\documentclass[leqno,a4paper,12pt]{article}
\usepackage{amsmath,amsthm,bbm}
\usepackage[sans]{dsfont}

\newtheorem{Thm}{\indent Theorem}[section]
\newtheorem{Prop}[Thm]{\indent Proposition}

\newtheorem{Cor}[Thm]{\indent Corollary}
\newtheorem{Var}[Thm]{\indent Variant}

\theoremstyle{definition}
\newtheorem{Def}[Thm]{\indent Definition}
\newtheorem{Rem}[Thm]{\indent Remark}
\newtheorem{Ex}[Thm]{\indent Example}
\newtheorem{Exo}[Thm]{\indent Exercice}
\newtheorem{Exs}[Thm]{\indent Examples}
\newtheorem{Cons}[Thm]{\indent Construction}

\newtheorem{Ques}[Thm]{\indent Question}
\def\qed{{\hskip0pt\unskip\unskip\nobreak\hfil\penalty50
          \hskip1em\hbox{}\nobreak\hfil
          {\bf q.e.d.}%
          \parfillskip=0pt\finalhyphendemerits=0
          \par}\medskip}

\newenvironment{Proof}
               {{\it Proof.}\quad}
               {\qed}

\newenvironment{Proofof}[1]
               {{\it Proof of #1.}\quad}
               {\qed}

\newcommand{\Prime}{\kern3\fontdimen1\font$'$\kern-7\fontdimen1\font}

\long\def\forget#1{}

\long\def\beginSIDEREMARK#1\endSIDEREMARK
    {{\par\bigskip\advance\leftskip by 2cm
                  \advance\rightskip by -2cm\noindent
      {\bf Our own side remark:} #1
      \par\bigskip\noindent}}
\long\def\beginFORGET#1\endFORGET{#1}
\long\def\beginFORGET#1\endFORGET{}

\def\?{\ ???\ \immediate\write16{}%
\immediate\write16{Warning: There was still a question mark . . . }%
\immediate\write16{}}

\usepackage{amsmath}
\usepackage{exscale}
\usepackage{amssymb}

\input cyracc.def
\font\tencyr=wncyr6
\def\cyr{\tencyr\cyracc}
\newcommand{\cyrb}{{\cyr B}}

\newcommand{\BA}{{\mathbb{A}}}

\newcommand{\BC}{{\mathbb{C}}}

\newcommand{\BG}{{\mathbb{G}}}

\newcommand{\BP}{{\mathbb{P}}}
\newcommand{\BQ}{{\mathbb{Q}}}
\newcommand{\BR}{{\mathbb{R}}}
\newcommand{\BS}{{\mathbb{S}}}

\newcommand{\BZ}{{\mathbb{Z}}}

\newcommand{\Fc}{{\mathfrak{c}}}

\newcommand{\Fn}{{\mathfrak{n}}}

\newcommand{\FA}{{\mathfrak{A}}}

\newcommand{\CC}{{\cal C}}
\newcommand{\CD}{{\cal D}}

\newcommand{\CH}{{\cal H}}

\newcommand{\CL}{{\cal L}}

\newcommand{\Spec}{\mathop{{\bf Spec}}\nolimits}

\newcommand{\rad}{\mathop{{\rm rad}}\nolimits}

\newcommand{\imm}{\mathop{{\rm im}}\nolimits}

\newcommand{\End}{\mathop{\rm End}\nolimits}

\newcommand{\GL}{\mathop{\rm GL}\nolimits}

\newcommand{\Gr}{\mathop{\rm Gr}\nolimits}

\newcommand{\Hom}{\mathop{\rm Hom}\nolimits}

\newcommand{\coker}{\mathop{\rm coker}\nolimits}

\newcommand{\loccit}{[loc.$\;$cit.]}

\def\halb{\frac{1}{2}}

\def\id{{\rm id}}

\newbox\mybox
\def\arrover#1{\mathrel{
       \setbox\mybox=\hbox spread 1.4em{\hfil$\scriptstyle#1$\hfil}
       \vbox{\offinterlineskip\copy\mybox
             \hbox to\wd\mybox{\rightarrowfill}}}}
\def\larrover#1{\mathrel{
       \setbox\mybox=\hbox spread 1.4em{\hfil$\scriptstyle#1$\hfil}
       \vbox{\offinterlineskip\copy\mybox
             \hbox to\wd\mybox{\leftarrowfill}}}}

\def\ontoover#1{\mathrel{
       \setbox\mybox=\hbox spread 1.4em{\hfil$\scriptstyle#1$\hfil}
       \vbox{\offinterlineskip\copy\mybox
             \hbox to\wd\mybox{\rightarrowfill\hskip-2.8mm
                               $\rightarrow$}}}}
\def\leftontoover#1{\mathrel{
       \setbox\mybox=\hbox spread 1.4em{\hfil$\scriptstyle#1$\hfil}
       \vbox{\offinterlineskip\copy\mybox
             \hbox to\wd\mybox{$\leftarrow$\hskip-2.8mm
                               \leftarrowfill}}}}
\def\longto{\longrightarrow}
\def\into{\hookrightarrow}
\def\onto{\ontoover{\ }}
\def\longonto{\ontoover{\ }}
\def\isoto{\arrover{\sim}}

\def\longinto{\lhook\joinrel\longrightarrow}

\usepackage[curve,matrix,arrow,cmtip]{xy}
\NoComputerModernTips

\def\myxymessage{\def\messagetext
   {Here an xy-pic diagram was omitted to speed up compilation . . . }
   \immediate\write16{\messagetext}
   \hbox{\bf \messagetext}}
\def\filxymatrix#1{\myxymessage}
\def\filxyarray#1{\myxymessage}
\newdir^{ (}{{}*!/-3pt/\dir^{(}}
\newdir_{ (}{{}*!/-3pt/\dir_{(}}
\newdir^{ )}{{}*!/+3pt/\dir^{)}}
\newdir_{ )}{{}*!/+3pt/\dir_{)}}

\def\rscript#1{\hbox to 0pt{$\scriptstyle#1$\hss}}

\let\oldbullet\bullet
\def\bullet{{\mathchoice{\oldbullet}%
                        {\oldbullet}%
                        {\scriptscriptstyle\oldbullet}%
                        {\oldbullet}}}

\newcommand{\argdot}{{\;\bullet\;}}%Punkt als Platzhalter fuer Argumente
\newcommand{\ujast}{\mathop{j_{!*}}\nolimits}

\newcommand{\CHM}{\mathop{CHM(k)}\nolimits}

\newcommand{\CHFAbM}{\mathop{CHM^{Ab}(k)_F}\nolimits}

\newcommand{\DBcM}{\mathop{DM_{\text{\cyrb},c}}\nolimits}

\newcommand{\DBcFkM}{\mathop{\DBcM(\Spec k)_F}\nolimits}

\newcommand{\DgM}{\mathop{DM_{gm}(k)}\nolimits}

\newcommand{\DAbgM}{\mathop{DM_{gm}^{Ab}(k)}\nolimits}
\newcommand{\DFAbgM}{\mathop{DM_{gm}^{Ab}(k)_F}\nolimits}

\newcommand{\MHM}{\mathop{\bf MHM}\nolimits}
\newcommand{\MHS}{\mathop{\bf MHS}\nolimits}

\newcommand{\Mgm}{\mathop{M_{gm}}\nolimits}

\newcommand{\Mcgm}{\mathop{M_{gm}^c}\nolimits}

\newcommand{\dm}{\mathop{\partial M}\nolimits}
\newcommand{\dMgm}{\mathop{\partial M_{gm}}\nolimits}

\newcommand{\one}{\mathds{1}}

\begin{document}

%%%%%%%%%%%%%%%%%%%%%%%%%%%%%%%%%%%%%%%%%%%%%%%%%%%%%%%%%%%%%%%%%%%%%%%
%
%  formatting

\hfuzz=3pt
\overfullrule=10pt                   % erzeugt schwarze Fehlerbalken

% The displayskip values were changed because \LaTeX does not react
% correctly to a \leqno: it should then use big skips, but doesn't.

\setlength{\abovedisplayskip}{6.0pt plus 3.0pt}
                               % preset 10.0pt plus 2.0pt minus 5.0pt
\setlength{\belowdisplayskip}{6.0pt plus 3.0pt}
                               % preset 10.0pt plus 2.0pt minus 5.0pt
\setlength{\abovedisplayshortskip}{6.0pt plus 3.0pt}
                               % preset 0.0pt plus 3.0pt
\setlength{\belowdisplayshortskip}{6.0pt plus 3.0pt}
                               % preset 6.0pt plus 3.0pt minus 3.0pt

\setlength{\baselineskip}{13.0pt}
                               % preset 12.0pt
\setlength{\lineskip}{0.0pt}
                               % preset 1.0pt
\setlength{\lineskiplimit}{0.0pt}
                               % preset 0.0pt

%%%%%%%%%%%%%%%%%%%%%%%%%%%%%%%%%%%%%%%%%%%%%%%%%%%%%%%%%%%%%%%%%%%%%%%
%
%  Title Page
%
%%%%%%%%%%%%%%%%%%%%%%%%%%%%%%%%%%%%%%%%%%%%%%%%%%%%%%%%%%%%%%%%%%%%%%%

\title{Absolute intersection motive
\forget{
\footnotemark
\footnotetext{To appear in ....}
}
}
\author{\footnotesize by\\ \\
\mbox{\hskip-2cm
\begin{minipage}{6cm} \begin{center} \begin{tabular}{c}
J\"org Wildeshaus \\[0.2cm]
\footnotesize Universit\'e Paris 13\\[-3pt]
\footnotesize Sorbonne Paris Cit\'e \\[-3pt]
\footnotesize LAGA, CNRS (UMR~7539)\\[-3pt]
\footnotesize F-93430 Villetaneuse\\[-3pt]
\footnotesize France\\
{\footnotesize \tt wildesh@math.univ-paris13.fr}
\end{tabular} \end{center} \end{minipage}
\hskip-2cm}
\\[2.5cm]
%{\bf Preliminary version --- not for distribution!}\\[1cm]
}
% In the final version we might want to fix the date:
%\date{March 18, 2011}
\maketitle
%\quad \\[-1.7cm]
\begin{abstract}
\noindent
The purpose of this article is
to define and study the notion of absolute intersection motive.  \\


%\noindent
%{\bf R\'esum\'e~:} RESUME.\\
\end{abstract}

%\vfill

\bigskip
\bigskip
\bigskip

\noindent {\footnotesize Math.\ Subj.\ Class.\ (2010) numbers: 14F32
(14C15, 14C25, 14C30, 14F25, 19E15, 32S60).
}

\eject

\tableofcontents

\bigskip
%\vspace*{0.5cm}

%\newpage

%\include{Intro}
%%%%%%%%%%%%%%%%%%%%%%%%%%%%%%%%%%%%%%%%%%%%%%%%%%%%%%%%%%%%%%%%%%%%%%%
%
%  Introduction
%
%%%%%%%%%%%%%%%%%%%%%%%%%%%%%%%%%%%%%%%%%%%%%%%%%%%%%%%%%%%%%%%%%%%%%%%

\setcounter{section}{-1}
\section{Introduction}
\label{Intro}

%%%%%%%%%%%%%%%%%%%%%%%%%%%%%%%%

%%%%%%%%%%%%%%%%%%%%%%%%%%%%%%%%

Let $X$ be a smooth scheme over the field $\BC$ of complex numbers.
According to Deligne, singular cohomology $H^n(X(\BC),\BQ)$ carries a \emph{mixed Hodge structure}
of weights $\ge n$, for all $n \in \BZ$. Dually, cohomology with compact support 
$H^n_c(X(\BC),\BQ)$ is equipped with a mixed Hodge structure of weights $\le n$. The canonical morphism
\[
H^n_c(X(\BC),\BQ) \longto H^n(X(\BC),\BQ)
\] 
is a morphism of Hodge structures. It therefore factors over a morphism 
\[
u_n: \Gr_n^W H^n_c(X(\BC),\BQ) \longto \Gr_n^W H^n(X(\BC),\BQ)
\]
between \emph{pure Hodge structures} of weight $n$,
where $\Gr_n^W H^n_c(X(\BC),\BQ)$ denotes the quotient of $H^n_c(X(\BC),\BQ)$ by
its part of weights $\le n-1$, and $\Gr_n^W H^n(X(\BC),\BQ)$ denotes  the part of $H^n(X(\BC),\BQ)$
of weight $n$. Therefore, the direct sum 
\[ 
\CH(u_n) := \ker(u_n) \oplus \imm(u_n) \oplus \coker(u_n) 
\]
is again a pure Hodge structure of weight $n$. \\

If $X$ is proper, then $H^n(X(\BC),\BQ) = H^n_c(X(\BC),\BQ)$ is pure, $u_n$ is an isomorphism,
and $\CH(u_n) = H^n(X(\BC),\BQ)$. Furthermore, 
the collection of all $\CH(u_n)$, $n \in \BZ$ is \emph{of motivic origin}:
there is a \emph{Chow motive} $\Mgm(X)$ over $\BC$, the \emph{motive of $X$},
whose (cohomological) \emph{Hodge theoretic realization} equals $(\CH(u_n))_{n \in \BZ}$. \\

The aim of this paper is to provide evidence for the following: whenever $X$ is smooth
(but not necessarily proper), then $(\CH(u_n))_{n \in \BZ}$ is of motivic origin: there is 
a Chow motive $M^{!*}(X)$, 
the \emph{absolute intersection motive of $X$}, 
whose Hodge theoretic realization equals $(\CH(u_n))_{n \in \BZ} \,$. \\

The $(\CH(u_n))_{n \in \BZ}$ satisfy the following minimality property: 
the morphism $u_n$ can be represented
as the composition of a monomorphism and an epimorphism
\[
\Gr_n^W H^n_c(X(\BC),\BQ) \longinto \CH(u_n) \longonto \Gr_n^W H^n(X(\BC),\BQ) \; ,
\]
and whenever $u_n$ is factored through a pure polarizable Hodge structure $H$ of weight $n$ in such a way,
then $\CH(u_n)$ is a direct factor of $H$. \\

The definition of the absolute intersection motive
$M^{!*}(X)$ (Definition~\ref{2C})
is conditional, but makes sense for any perfect base field $k$.
Indeed, it depends on the existence of a \emph{minimal weight filtration} of the \emph{boundary motive}
$\dMgm(X)$. This existence would be ensured if as expected, all Chow motives over $k$,
or at least those necessary to build up $\dMgm(X)$, 
were \emph{finite dimensional} in the sense of Kimura. A priori, given the main result from \cite{Ki}, this
means that $M^{!*}(X)$ can be defined as soon as $\dMgm(X)$ is \emph{of Abelian type}. \\

Section~\ref{1} contains the necessary basics from the theory of minimal weight filtrations.
With an eye to the application to the Hodge theoretic context, we take particular care to formulate
results for triangulated categories, which are not only equipped with a \emph{weight structure} $w$
\`a la Bondarko, but equally with a $t$-structure \emph{transversal} to $w$. The main result is Theorem~\ref{1F},
which gives a characterization of morphisms in the \emph{radical} in terms of their effect on cohomology objects. \\ 

Section~\ref{2} contains the basic definitions, conditional as we said in the motivic context, but unconditional
in the Hodge theoretic one. We thus get the notion of \emph{absolute Hodge theoretic intersection complex}
of a smooth $\BC$-scheme $X$ (Definition~\ref{2D}). Both definitions rely on
Construction~\ref{Cons}, which establishes a bijection between isomorphism classes 
of \emph{weight filtrations} of some boundary object on the one hand,
and isomorphism classes of fac\-torizations. \\

By definition,
the cohomology objects of the absolute Hodge theoretic intersection complex give \emph{absolute intersection
cohomology}; its study is the object of Section~\ref{3}. In particular (Theorem~\ref{3H}), we prove,
using the theo\-ry developed in Section~\ref{1}, that
absolute intersection cohomology of $X$
is isomorphic to $(\CH(u_n))_{n \in \BZ} \,$; therefore it satisfies the minimality property
mentioned above (Corollary~\ref{3K}). \\

Section~\ref{4} contains a great number of examples. We establish
the first examples of minimal weight filtrations for motives which are
not of Abelian type (Example~\ref{4E}, Corollary~\ref{4Fa}). We also see that
in general, the absolute intersection motive
of a product is unequal to the tensor product of the absolute intersection motives
of the factors (Remark~\ref{4Ga}). \\

Section~\ref{5} is concerned with the following problem (Question~\ref{5A}): can absolute intersection cohomology
be ``geometrically realized'', \emph{i.e.}, does there exist a topological space $X^{!*}$, stratified into
topological manifolds, among which $X(\BC)$, such that absolute intersection cohomology of $X$ equals
intersection cohomology relative to $X^{!*} \,$? None of the examples we found suggests the contrary, and
we identify $X^{!*}$ in many cases. Actually, we find that sometimes, but not always, the Alexandrov one-point
compactification $X^+$ can be chosen as $X^{!*}$. Theorem~\ref{5E} gives
necessary and sufficient criteria for the equation ``$X^{!*} = X^+$'' to hold. They basically tell us that
whenever $Y$ is a smooth algebraic compactification of $X$, with complement $Z$, then ``$X^{!*} = X^+$''
if and only if ``the complement $Z$ has maximal geometric self-interaction''. By contrast, ``no geometric
self-interaction of $Z$'' seems to lead to the solution ``$X^{!*} = Y(\BC)$''. We finish the article
with an example (Exercice~\ref{5F}), where $X^{!*}$ exists, but is neither equal to $X^+$
nor to $Y(\BC)$, for any smooth compactification $Y$ of $X$. \\  

Part of this work was done while I was enjoying a \emph{d\'el\'egation
aupr\`es du CNRS}, and during a visit to Institut Mittag-Leffler (Djursholm). 
I am grateful to both institutions. \\

{\bf Notation and conventions}: For a perfect field $k$, 
we denote by $Sch/k$ the category of separated schemes of finite 
type over $k$, and by $Sm/k \subset Sch/k$
the full sub-category of objects which are smooth over $k$. 
As far as motives are concerned,
the notation of this paper is that of \cite{W2,W4}, which in turn follows
that of \cite{V}. We refer to \cite[Sect.~1]{W2} for a concise
review of this notation, and of the de\-fi\-nition of the triangulated 
category $\DgM$ of geometrical
motives over $k$. Let $F$ be a commutative semi-simple Noetherian $\BQ$-algebra,
in other words, a finite direct product of fields of characteristic zero.
The notation $\DgM_F$ stands for the 
$F$-linear analogue of $\DgM$
defined in \cite[Sect.~16.2.4
and Sect.~17.1.3]{A}. 
Simi\-larly,
let us denote by $\CHM$ the category opposite to the category
of Chow motives, and by
$\CHM_F$ the pseudo-Abelian
completion of the category
$\CHM \otimes_\BZ F$. Using \cite[Cor.~2]{V2}, we canonically identify $\CHM_F$ with
a full additive sub-category of $\DgM_F \,$. \\

When we assume a field $k$ to \emph{admit resolution of singularities},
then it will be 
in the sense of \cite[Def.~3.4]{FV}:
(i)~for any separated $k$-scheme $X$ of finite type, there exists an abstract blow-up $Y \to X$
\cite[Def.~3.1]{FV} whose source $Y$ is smooth over $k$,
(ii)~for any pair of smooth, seperated $k$-schemes $X, Y$ of finite type, 
and any abstract blow-up $q : Y \to X$,
there exists a sequence of blow-ups 
$p: X_n \to \ldots \to X_1 = X$ with smooth centers,
such that $p$ factors through $q$. 
We say that $k$ \emph{admits strict resolution of singularities},
if in (i), for any given dense open subset $U$
of the smooth locus of $X$,
the blow-up $q: Y \to X$ can be chosen to be an isomorphism above $U$,
and such that arbitrary intersections of
the irreducible components of the complement $Z$ of $U$ in $Y$  
are smooth (e.g., $Z \subset Y$ a normal
crossing divisor with smooth irreducible components).

%%% Local Variables:
%%% mode: latex
%%% TeX-master: "head"
%%% End:

\bigskip
%\include{Sec1}
%%%%%%%%%%%%%%%%%%%%%%%%%%%%%%%%%%%%%%%%%%%%%%%%%%%%%%%%%%%%%%%%%%%%%%%
%
%  Section 1
%
%%%%%%%%%%%%%%%%%%%%%%%%%%%%%%%%%%%%%%%%%%%%%%%%%%%%%%%%%%%%%%%%%%%%%%%

\section{Weight structures and $t$-structures}
\label{1}

%%%%%%%%%%%%%%%%%%%%%%%%%%%%%%%%

%%%%%%%%%%%%%%%%%%%%%%%%%%%%%%%%

We make free use of the terminology of and basic results on weight structures
\cite[Sect.~1.3]{B2}.  Let us fix an $F$-linear triangulated cate\-go\-ry $\CC$,
which is equipped with a weight structure $w = (\CC_{w \le 0},\CC_{w \ge 0})$.
For simplicity, and in order to be able to apply the results from \cite{B3},
we assume $w$ to be bounded.
The following notion is the key for everything to follow.

\begin{Def}[{\cite[Def.~1.3]{W11}}] \label{1A}
Let $M \in \CC$, and $n \in \BZ$. 
A \emph{minimal weight filtration concentrated at $n$} of $M$ 
is a weight filtration
\[
M_{\le n-1} \longto M \longto M_{\ge n} \stackrel{\delta}{\longto} M_{\le n-1}[1]
\]
($M_{\le n-1} \in \CC_{w \le n-1}$, $M_{\ge n} \in \CC_{w \ge n}$) 
such that the morphism $\delta$ belongs to the \emph{radical} 
\cite[D\'ef.~1.4.1]{AK} of $\CC$: 
\[
\delta \in \rad_\CC (M_{\ge n},M_{\le n-1}[1]) \; .
\]
\end{Def}

Any two minimal weight filtrations of the same object $M$ are related
by an isomorphism (which in general is \emph{not} unique) 
\cite[proof of Thm.~2.2~(b)]{W9}. Minimal weight filtrations do not necessarily exist
\cite[Ex.~2.3~(c)]{W9}. 

\begin{Rem} \label{1B}
(a)~Minimal weight filtrations \emph{do} exist if the \emph{heart} 
\[
\CC_{w = 0} := \CC_{w \le 0} \cap \CC_{w \ge 0}
\] 
of $w$ is pseudo-Abelian and \emph{semi-primary} \cite[Thm.~2.2~(a)]{W9}, 
\emph{i.e.}~\cite[D\'ef.~2.3.1]{AK}, if \\[0.1cm]
(1)~for all objects $M$ of $\CC_{w = 0}$, the radical $\rad_{\CC_{w = 0}} (M,M)$ 
is nilpotent, \\[0.1cm]
(2)~the $F$-linear quotient category $\CC_{w = 0} / \rad_{\CC_{w = 0}}$
is semi-simple. \\[0.1cm]
(b)~In practice, given $M \in \CC$, the existence of minimal weight
filtrations of $M$ is assured once $M$ can be shown to belong to a full,
triangulated sub-category $\CC'$ of $\CC$, such that $w$ induces a weight structure
on $\CC'$, and such that $\CC_{w = 0}' = \CC' \cap \CC_{w = 0}$ is pseudo-Abelian and semi-primary.
\end{Rem}

\begin{Ex} \label{1Ba}
According to  \cite[Sect.~6]{B1}, 
the category $\CC := \DgM_F$ of geometrical motives over $k$
carries a bounded weight structure $w$, if $k$ is of characteristic zero. This claim still holds for arbitrary
perfect fields (remember that $F$ is supposed to be a $\BQ$-algebra), as can be seen from 
the proof of \cite[Thm.~1.13]{W4}, using the main results from \cite[Sect.~5.5]{Ke}. 
We refer to this weight structure as \emph{motivic}.
It identifies $\CHM_F$ with the heart $\DgM_{F,w=0} \,$.  
The category $\CHM_F$ is pseudo-Abelian. It is expected, but not known to be semi-primary.
Consider the full sub-category $\CHFAbM$ of $\CHM_F$ of \emph{Chow motives of Abelian
type over $k$} \cite[Def.~1.1~(b)]{W9}, 
and the full triangulated sub-category $\CC':= \DFAbgM$ of $\DgM_F$
generated by $\CHFAbM \,$. 
The motivic weight structure induces a weight structure, still denoted $w$, on $\DFAbgM \,$
\cite[Prop.~1.2]{W9},
and $DM_{gm}^{Ab}(k)_{F,w=0} = \CHFAbM \,$. According to \cite[Prop.~1.8]{W9},
the category $\CHFAbM$ is (pseudo-Abelian and) semi-primary. 
We may thus apply Remark~\ref{1B}~(b): any geometrical motive belonging 
to $\DFAbgM$ admits minimal weight filtrations.
\end{Ex}

In the sequel of this section, we shall
assume in addition that $\CC$ carries a bounded $t$-structure 
$t = (\CC^{t \le 0},\CC^{t \ge 0})$, which is \emph{transversal}
to $w$ \cite[Def.~1.2.2]{B3}. (Note that in the motivic context discussed in Example~\ref{1Ba},
such a $t$-structure is expected, but not known to exist.)
Denote by $\CC_{w = n}^{t = m} \,$, $m, n \in \BZ$,
the full sub-categories given by the intersections
\[
\CC_{w = n}^{t = m} := \CC_{w = n} \cap \CC^{t = m}
\]
(where $\CC_{w = n} := \CC_{w \le n} \cap \CC_{w \ge n} = \CC_{w = 0}[n]$
and $\CC^{t = m} := \CC^{t \le m} \cap \CC^{t \ge m} = \CC^{t = 0}[-m]$).

\begin{Prop}[{\cite[Thm.~1.2.1, Rem.~1.2.3, 
Prop.~1.2.4]{B3}}] \label{1C}
Assume $\CC$ to be equally equipped with a bounded $t$-structure, 
which is transversal to $w$. \\[0.1cm]
(a)~The categories $\CC_{w = n}^{t = m} \,$, $m,n \in \BZ$, are Abelian semi-simple,
and any object of $\CC_{w = 0}$
is isomorphic to a finite direct sum of objects in $\CC_{w = 0}^{t = m} \,$,
$m \in \BZ$. \\[0.1cm]
(b)~For fixed $m \in \BZ$, and $n_1 \ne n_2$,
there are no non-zero morphisms between objects
of $\CC_{w = n_1}^{t = m}$ and of $\CC_{w = n_2}^{t = m}$. \\[0.1cm]
(c)~For fixed $m \in \BZ$, any object of $\CC^{t = m}$ admits weight filtrations 
by objects of $\CC^{t = m}$, which are exact sequences in $\CC^{t = m}$. \\[0.1cm]
(d)~The $t$-truncation functors $\tau^{t \le m}$, $\tau^{t \ge m}$, $m \in \BZ$,
respect the sub-categories $\CC_{w \ge n}$ and $\CC_{w \le n}$, $n \in \BZ$.
\end{Prop} 

\begin{Cor} \label{1Ca}
Assume $\CC$ to be equally equipped with a bounded $t$-struc\-ture, 
which is transversal to $w$. \\[0.1cm]
(a)~The heart $\CC_{w = 0}$ is pseudo-Abelian and semi-primary. \\[0.1cm]
(b)~Any object of $\CC$ admits mi\-ni\-mal weight filtrations.
\end{Cor}

\begin{Proof}
As recalled in Remark~\ref{1B}~(a), part~(b) is implied by (a).

The second claim of part~(a) follows from Proposition~\ref{1C}~(a), and from
\cite[proof of Prop.~2.3.4~c)]{AK}. It remains to show that $\CC_{w = 0}$ is pseudo-Abelian.
Let $M = \oplus_{m \in \BZ} M^m$ be an object of $\CC_{w = 0}$, with 
$M^m \in \CC_{w = 0}^{t = m} \,$, almost all $M^m$ being zero (Proposition~\ref{1C}~(a)).
Let $e$ be an idempotent endomorphism of $M$.
In order to show that $e$ admits a kernel, we apply induction on the number of 
non-zero components $M^m$, the initial case $M = M^m$ 
resulting from Proposition~\ref{1C}~(a). For the induction step, 
take $m$ to be minimal such that $M^m \ne 0$, and
write  
\[
M = M^m \oplus N \; .
\]  
Orthogonality for the $t$-structure tells us that with respect to this direct sum, 
\[
e =
\left( \begin{array}{cc}
A & B \\
0 & D
\end{array} \right) \; ,
\]
with $A \in \End_{\CC_{w = 0}} (M^m)$, $B \in \Hom_{\CC_{w = 0}} (N,M^m)$, and 
$D \in \End_{\CC_{w = 0}} (N)$.
The relation $e^2 = e$ is equivalent to the system of relations
\[
(\ast) \quad \quad A^2 = A \; , \; D^2 = D \; , \; AB + BD = B \; .
\] 
By our induction hypothesis, $\ker(A) \subset M^m$ and $\ker(D) \subset N$ exist
(and so do $\ker(\id_{M^m} - A)$ and $\ker(\id_N - D)$).
We leave it as an exercice to the reader to prove, using $(\ast)$, that
the morphism
\[
\left( \begin{array}{cc}
\id_{\ker(A)} & -B \\
0 & \id_{\ker(D)}
\end{array} \right) 
\]
from $\ker(A) \oplus \ker(D)$ to $M$ is a kernel of $e$.
\end{Proof}

Although it will not be needed in the sequel of the present article, let us
mention the following.

\begin{Cor}
Let $\CC$ be equipped with a bounded weight structure $w$, and 
with a bounded $t$-structure, which is transversal to $w$.
Then $\CC$ is pseudo-Abelian.
\end{Cor}

\begin{Proof}
This follows from Corollary~\ref{1Ca}~(a), and from \cite[Lemma~5.2.1]{B1}.
\end{Proof}

Write
\[
H^m : \CC \longto \CC^{t = 0} \; , \; m\in \BZ \; ,
\]
for the cohomology functors associated to $t$. 
The following holds even in the absence of a weight structure.

\begin{Prop} \label{1D}
Let $M$ and $N$ be objects of $\CC$, and $\alpha: M \to N$ a morphism.
If $H^m(\alpha) = 0$ for all $m \in \BZ$, then $\alpha \in \rad_\CC (M,N)$. 
\end{Prop}

\begin{Proof}
For any morphism $\beta: N \to M$, and any $m \in \BZ$, we have 
\[
H^m \bigl( \id_M - \beta \alpha \bigr) = \id_{H^m(M)} \; ,
\]
meaning that $H^*(\id_M - \beta \alpha)$ is an automorphism. Hence
so is $\id_M - \beta \alpha$.
\end{Proof}

\begin{Ex} \label{1E}
The converse of Proposition~\ref{1D} is not true in general. Let $N$ be a simple object 
of $\CC_{w = 0}^{t = 0} \,$ and $M$ a non-trivial extension in $\CC^{t = 0}$
\[
0 \longto N_- \longto M \stackrel{\alpha}{\longto} N \longto 0
\]
of $N$ by some object $N_-$ of $\CC^{t = 0}$ of strictly negative weights. 
Schur's Lemma and Proposition~\ref{1C}~(c), (b) imply that
\[
\Hom_\CC (N,M) = 0 \; .
\]
Therefore, $\alpha \in \rad_\CC (M,N)$, but $0 \ne \alpha = H^0 (\alpha)$. 
\end{Ex}

It turns out that the converse of Proposition~\ref{1D} \emph{is} true,
once we restrict the weights of $M$ and $N$. Here is the main result
of this section.

\begin{Thm} \label{1F}
Assume $\CC$ to be equally equipped with a bounded $t$-struc\-ture, 
which is transversal to $w$. Let $n \in \BZ$, $M \in \CC_{w \ge n}$,
and $N \in \CC_{w \le n}$. Let $\alpha: M \to N$ be a morphism.
Then the following are equivalent.
\begin{enumerate}
\item[(1)] $H^m(\alpha) = 0$ for all $m \in \BZ$.
\item[(2)] $\alpha \in \rad_\CC (M,N)$.
\end{enumerate} 
\end{Thm}

\begin{Proof}
Given Proposition~\ref{1D}, we may assume that $\alpha$ belongs to the radical,
and need to establish that it is zero on cohomology. 

Let us first treat the case where both $M$ and $N$ are objects of some $\CC^{t = r}$, and
show the following claim~$(\ast)$: if a morphism $\gamma: M \to N$ is non-zero,
then there is $\beta: N \to M$ such that $\beta \gamma: M \to M$ is a projector
onto a non-zero direct factor of $M$. 

Indeed, let
\[
0 \longto W_n M \longto M \longto W^{n+1} M \longto 0
\]
and 
\[
0 \longto W_{n-1} N \longto N \longto W^n N \longto 0
\]
be weight filtrations of $M$ and $N$ within $\CC^{t = r}$ (Proposition~\ref{1C}~(c)),
with $W_n M$ and $W^n N$ in $\CC_{w = n}^{t = r} \,$, $W^{n+1} M \in \CC_{w \ge n+1}$,
and $W_{n-1} N \in \CC_{w \le n-1}$. Write $\gamma_n$ for the composition
\[
W_n M \longinto M \stackrel{\gamma}{\longto} N \longonto W^n N \; ;
\]
Proposition~\ref{1C}~(c), (b) imply that $\gamma_n \ne 0$.
According to Proposition~\ref{1C}~(a), there is $\beta_n : W^n N \to W_n M$
such that $\beta_n \gamma_n$ is a projector onto a complement of 
$\ker (\gamma_n) \subset W_n M$. 
Define $\beta$ as the composition
\[
N \longonto W^n N \stackrel{\beta_n}{\longto} W_n M \longinto M \; .
\]
Now let us come back to $\alpha \in \rad_\CC (M,N)$. In order to show that
$H^m(\alpha) = 0$ for all $m \in \BZ$, let us apply induction on the sum $s$ of
the number of non-zero cohomology objects of $M$ and of $N$. The claim is trivial if
$s \le 1$, or if there is no degree $r$ for which both $H^r(M)$ and $H^r(N)$
are non-zero. In order to treat the case $s = 2$, we may thus assume that
$M, N \in \CC^{t = r}$, for some $r \in \BZ$. Claim~$(\ast)$ then implies that indeed
\[
0 = \alpha = H^r(\alpha) \; .
\]
For the induction step, let $r \in \BZ$ be maximal such that $H^r(M)$ and $H^r(N)$ 
are both non-zero. This takes trivially care of $H^m(\alpha)$, for all $m \ge r+1$.

Using $t$-truncations $\tau^{t \le \argdot}$, $\tau^{t \ge \argdot}$, we see 
that the induction hypothesis can be applied unless 
$M \in \CC^{t \le r}$ and $N \in \CC^{t \ge r}$. In that case, 
$H^m(\alpha) = 0$ for all $m \ne r$, and 
$\alpha$ equals the composition
\[
M \longto H^r(M)[-r] \stackrel{\gamma}{\longto} H^r(N)[-r] \longto N \; ,
\]
with $\gamma := H^r(\alpha)[-r]$. The objects $H^r(M)[-r]$ and $H^r(N)[-r]$
belong to $\CC_{w \ge n}$ and $\CC_{w \le n}$, respectively 
(Proposition~\ref{1C}~(d)). If $\gamma$ were non-zero, then claim~$(\ast)$ would show that
there is $\beta': H^r(N)[-r] \to H^r(M)[-r]$ such that $\beta' \gamma$ is a
projector onto a non-trivial direct factor of $H^r(M)[-r]$. The maps
\[
\Hom_\CC (N,M) \longto \Hom_\CC (H^r(N)[-r],M)
\]
and
\[
\Hom_\CC (H^r(N)[-r],M) \longto \Hom_\CC (H^r(N)[-r],H^r(M)[-r])
\]
are both surjective, since $\tau^{t \ge r+1} (N)[-1] \in \CC_{w \le n-1}$
and $\tau^{t \le r-1} (M)[1] \in \CC_{w \ge n+1}$ (Proposition~\ref{1C}~(d)). 
Thus, the morphism $\beta'$ can be extended to yield $\beta: N \to M$.
By construction, the composition $\beta \alpha: M \to M$ is a projector
onto a non-trivial direct factor of $M$, implying in particular that
$\id_M - \beta \alpha$ is not an automorphism of $M$. But this contradicts
the assumption on $\alpha$.
\end{Proof}
 
Theorem~\ref{1F} has the following important consequences.

\begin{Cor} \label{1G}
Assume $\CC$ to be equally equipped with a bounded $t$-struc\-ture, 
which is transversal to $w$. Let 
\[
M_{\le n-1} \longto M \longto M_{\ge n} \stackrel{\delta}{\longto} M_{\le n-1}[1]
\]
be a minimal weight filtration concentrated at some $n \in \BZ$ of an object $M$ of $\CC$
(according to Corollary~\ref{1Ca}(b), such minimal weight filtrations exist).
Then for all $m \in \BZ$, the sequence 
\[
0 \longto H^m \bigl( M_{\le n-1} \bigr) \longto H^m(M) 
\longto H^m \bigl( M_{\ge n} \bigr) \longto 0
\]
is exact, and yields the weight filtration of $H^m(M)$ concentrated at $m+n$.
In particular, the sequence does not depend on the choice of minimal weight filtration.
\end{Cor}

\begin{Proof}
Apply Theorem~\ref{1F} to $\delta: M_{\ge n} \to M_{\le n-1}[1]$.
Thus, the long exact cohomology sequence yields short exact sequences 
\[
0 \longto H^m \bigl( M_{\le n-1} \bigr) \longto H^m(M) 
\longto H^m \bigl( M_{\ge n} \bigr) \longto 0 \; ,
\]
$m \in \BZ$. According to Proposition~\ref{1C}~(d), $H^m ( M_{\le n-1} )[-m]$
is indeed of weights $\le n-1$, and $H^m ( M_{\ge n} )[-m]$, of weights $\ge n$. 
\end{Proof}

\begin{Cor} \label{1H}
Assume $\CC$ to be equally equipped with a bounded $t$-struc\-ture, 
which is transversal to $w$. Let $M, N \in \CC_{w = 0}$.
Then
\[
\rad_{\CC_{w = 0}} (M,N) = \{ \alpha \in \Hom_{\CC_{w = 0}} (M,N) \; , \; 
H^m(\alpha) = 0 \;\; \text{for all} \;\; m \in \BZ \} \; .
\]
\end{Cor}

\begin{Ex} \label{1I}
Assume $F$ to be a field contained
in $\BR$, and put
$\CC := D^b ( \MHS_F )$, the bounded derived category 
of \emph{mixed graded-polarizable $F$-Hodge structures}
\cite[Def.~3.9, Lemma~3.11]{Be}. 
According to \cite[Prop.~2.3.1~I]{B3} (with $X = \Spec \BC$), $\CC$ carries a canonical
weight structure $w$, which is bounded: indeed, the 
class of a bounded complex $K$ of mixed graded-pola\-rizable $F$-Hodge structures
lies in $D^b ( \MHS_F )_{w \le 0}$ (resp., in $D^b ( \MHS_F )_{w \ge 0}$) if and
only if the $m$-th cohomo\-lo\-gy object $H^m(K)$ is a Hodge structure of weights
$\le m$ (resp., $\ge m$), for all $m \in \BZ$.
Furthermore, the canonical $t$-structure
on $D^b ( \MHS_F )$ is transversal to $w$. Therefore, the theory developed 
in the present section applies. In particular
(Corollary~\ref{1Ca}~(a)), $D^b ( \MHM_F X )_{w = 0}$ 
is pseudo-Abelian and semi-primary. 
\end{Ex}

\begin{Ex} \label{1K}
Examples~\ref{1Ba} and \ref{1I} are related by the \emph{Hodge theoretic realization}
(\cite[Sect.~2.3 and Corrigendum]{Hu}; see \cite[Sect.~1.5]{DG}
for a simplification of this approach). The field $k$ is assumed to be
embedded into $\BC$ \emph{via} $\sigma: k \into \BC$, and 
$F$ is assumed to be a sub-field of $\BR$. 
Denote by 
\[
R_\sigma : \DgM_F \longto D^b ( \MHS_F )
\]
the Hodge theoretic realization of \loccit , and
recall that it is a \emph{contrava\-riant} tensor functor mapping
the pure Tate motive $\BZ(m)$ to the pure Hodge structure $\BQ(-m)$
\cite[Thm.~2.3.3]{Hu}. Chow motives are mapped to objects of $D^b ( \MHS_F )$,
which are pure of weight zero. Since the motivic weight structure is bounded,
it follows formally that 
\[
R_\sigma \bigl( \DgM_{F,w \le 0} \bigr) \subset D^b ( \MHS_F )_{w \ge 0} \; ,
\]
and that
\[
R_\sigma \bigl( \DgM_{F,w \ge 0} \bigr) \subset D^b ( \MHS_F )_{w \le 0} \; .
\]
\end{Ex}

In the setting of Example~\ref{1K}, we have the following result.

\begin{Thm} \label{1L}
Consider the restriction 
\[
R_\sigma^{Ab} : \DFAbgM \longto D^b ( \MHS_F )
\]
of $R_\sigma$ to $\DFAbgM \subset \DgM_F$. \\[0.1cm]
(a)~Let $n \in \BZ$, $M \in \DAbgM_{F,w \ge n}$,
and $N \in \DAbgM_{F,w \le n}$. Then
\[
R_\sigma^{Ab} \bigl( \rad_{\DFAbgM} (M,N) \bigr) \subset 
\rad_{D^b ( \MHS_F )} \bigl(( R_\sigma(N) , R_\sigma(M) \bigr) \; .
\]
(b)~The functor $R_\sigma^{Ab}$ maps minimal weight filtrations (concentrated at $n$)
to minimal weight filtrations (concentrated at $-n+1$).
\end{Thm}

\begin{Proof}
Part~(a) follows from \cite[Lemma~1.11]{W8} and Proposition~\ref{1D}.

Given (a), part~(b) follows from the definition (and the contravariance of $R_\sigma^{Ab}$).
\end{Proof} 

\begin{Rem}
The importance of the hypothesis ``of Abelian type'' is twofold. 
First, it guarantees the existence of minimal weight filtrations of objects of
$\DFAbgM$. Note that this would still be true if we replaced $\DFAbgM$ by the full
triangulated sub-category of $\DgM_F$ generated by Chow motives, which 
are \emph{finite dimensional} \cite[Def.~3.7]{Ki}. 
Second, and more seriously, Theorem~\ref{1L}~(a) implies that the restriction
of $R_\sigma$ to $\CHFAbM$ is \emph{radicial} \cite[D\'ef.~1.4.6]{AK}. Indeed,
this latter statement is the vital ingredient of the proof of 
Theorem~\ref{1L}~(a); see \cite[Thm.~1.9]{W8}. But as shown by the proof of \loccit,
that statement is intimately related to the identification of homological and numerical
equivalence, which thanks to \cite{L} 
is known for motives associated to Abelian varieties. 
\end{Rem}

%%% Local Variables:
%%% mode: latex
%%% TeX-master: "head"
%%% End:

\bigskip
%\include{Sec2}
%%%%%%%%%%%%%%%%%%%%%%%%%%%%%%%%%%%%%%%%%%%%%%%%%%%%%%%%%%%%%%%%%%%%%%%
%
%  Section 2
%
%%%%%%%%%%%%%%%%%%%%%%%%%%%%%%%%%%%%%%%%%%%%%%%%%%%%%%%%%%%%%%%%%%%%%%%

\section{Definitions}
\label{2}

%%%%%%%%%%%%%%%%%%%%%%%%%%%%%%%%

%%%%%%%%%%%%%%%%%%%%%%%%%%%%%%%%

Let us come back to the situation from the beginning of Section~\ref{1},
\emph{i.e.}, an $F$-linear triangulated cate\-go\-ry $\CC$,
equipped with a bounded weight structure $w = (\CC_{w \le 0},\CC_{w \ge 0})$
(but not necessarily with a $t$-structure). 
Let us fix a morphism $u: M_- \to M_+$ in $\CC$ between 
$M_- \in \CC_{w \le 0}$ and $M_+ \in \CC_{w \ge 0} \ $.

\begin{Cons} \label{Cons}
First choose a cone $C$ of $u$, \emph{i.e.}, an exact triangle
\[
C[-1] \longto M_- \stackrel{u}{\longto}
M_+ \longto C 
\]
in $\CC$. 
Then choose a weight filtration of $C$
\[
C_{\le 0} \longto C \longto C_{\ge 1} 
\stackrel{\delta}{\longto} C_{\le 0}[1]
\]
(with $C_{\le 0} \in \CC_{w \le 0}$ and $C_{\ge 1} \in \CC_{w \ge 1}$).
Consider the diagram, which we shall denote by the symbol~(1), of exact triangles 
\[
\vcenter{\xymatrix@R-10pt{
        0 \ar[d] \ar[r] &
        C_{\ge 1}[-1] \ar@{=}[r] &
        C_{\ge 1}[-1] \ar[d]^{\delta[-1]} \ar[r] &
        0 \ar[d] \\
        M_- \ar@{=}[d] &  
         &
        C_{\le 0} \ar[d] \ar[r] &
        M_-[1] \ar@{=}[d] \\
        M_- \ar[d] \ar[r]^-u &
        M_+ \ar[d] \ar[r] &
        C \ar[d] \ar[r] &
        M_-[1] \ar[d] \\
        0 \ar[r] &
        C_{\ge 1} \ar@{=}[r] &
        C_{\ge 1} \ar[r] &
        0    
\\}}
\]
which according to axiom TR4' of triangulated categories 
\cite[Sect.~1.1.6]{BBD}
can be completed to give a diagram, denoted by~(2)
\[
\vcenter{\xymatrix@R-10pt{
        0 \ar[d] \ar[r] &
        C_{\ge 1}[-1] \ar[d]_{\delta_+[-1]} \ar@{=}[r] &
        C_{\ge 1}[-1] \ar[d]^{\delta[-1]} \ar[r] &
        0 \ar[d] \\
        M_- \ar@{=}[d] \ar[r]^-i &  
        M \ar[d]_{\pi} \ar[r]^-{\delta_-[-1]} &
        C_{\le 0} \ar[d] \ar[r] &
        M_-[1] \ar@{=}[d] \\
        M_- \ar[d] \ar[r]^-u &
        M_+ \ar[d] \ar[r] &
        C \ar[d] \ar[r] &
        M_-[1] \ar[d] \\
        0 \ar[r] &
        C_{\ge 1} \ar@{=}[r] &
        C_{\ge 1} \ar[r] &
        0    
\\}}
\]
with $M \in \CC$. By the second row, and the second column of
diagram~(2), the object $M$
is simultaneously an extension of objects of weights $\le 0$, 
and an extension of objects of weights $\ge 0$.
It follows easily (\emph{c.f.}\ \cite[Prop.~1.3.3~3]{B1})
that $M$ belongs to both $\CC_{w \le 0}$ and $\CC_{w \ge 0} \, $, and hence to
$\CC_{w=0} \, $.

Note that even for a fixed choice of $C$, 
diagram~(2) necessitates further choices of $M$
and of factorizations $u = \pi i$ and $\delta = \delta_- \delta_+$.
In general, the object $M$ is unique up to possibly non-unique
isomorphism.
\end{Cons}
 
Very precisely, we have the following.

\begin{Prop} \label{2A}
The map
\[
\big\{ \bigl( C_{\le 0} , C_{\ge 1} \bigr) \big\} / \cong \; 
\longto 
\big\{ \bigl( M , i , \pi \bigr) \big\} / \cong
\]
from the Construction~\ref{Cons}
is a bijection between \\[0.1cm]
(1)~the isomorphism classes of cones $C$ of $u$, together with weight filtrations of $C$, \\[0.1cm]
(2)~the isomorphism classes of objects
$M$ of $\CC_{w=0} \, $, together with a factorization 
\[
M_- \stackrel{i}{\longto} M
\stackrel{\pi}{\longto} M_+
\]
of the morphism $u: M_- \to M_+$, such that 
both $i$ and $\pi$ can be completed to give weight filtrations of 
$M_-$ and of $M_+$, respectively.
\end{Prop}

\begin{Def} \label{2B}
Assume that a cone $C$ of $u: M_- \to M_+$ (hence \emph{any} cone of $u$)
admits a minimal weight filtration 
\[
C_{\le 0} \longto C \longto C_{\ge 1} 
\stackrel{\delta}{\longto} C_{\le 0}[1]
\]
concentrated at $n=1$. A \emph{minimal factorization of $u: M_- \to M_+$} 
is a triple $( M , i , \pi )$ consisting 
of an object
$M$ of $\CC_{w=0} \, $, together with a factorization 
\[
M_- \stackrel{i}{\longto} M
\stackrel{\pi}{\longto} M_+
\]
of $u$, whose isomorphism class corresponds to $(C_{\le 0} , C_{\ge 1})$
under the bijection of Proposition~\ref{2A}.
\end{Def}

Thus, given two minimal factorizations 
\[
M_- \stackrel{i_1}{\longto} M_1 \stackrel{\pi_1}{\longto} M_+ 
\quad \text{and} \quad  
M_- \stackrel{i_2}{\longto} M_2 \stackrel{\pi_2}{\longto} M_+ 
\]
of $u$, there is an isomorphism (which in general is not unique) $\alpha: M_1 \isoto M_2$
such that the diagram
\[
\vcenter{\xymatrix@R-10pt{
        M_- \ar@{=}[d] \ar[r]^-{i_1} &  
        M_1 \ar[d]_\alpha \ar[r]^-{\pi_1} &
        M_+ \ar@{=}[d] \\
        M_- \ar[r]^-{i_2} &  
        M_2 \ar[r]^-{\pi_2} &
        M_+
\\}}
\]
commutes. 

\begin{Rem} \label{2Ba}
(a)~The correspondence from Proposition~\ref{2A} is not \emph{a priori} functorial,
meaning that a morphism between weight filtrations of $C$ does not yield a morphism 
of the corresponding objects $M$, not even up to automorphism. \\[0.1cm]
(b)~Minimality of the weight filtration
\[
C_{\le 0} \longto C \longto C_{\ge 1} 
\stackrel{\delta}{\longto} C_{\le 0}[1]
\]
of $C$ does in general not imply minimality of the weight filtrations
\[
C_{\le 0}[-1] \longto M_- \stackrel{i}{\longto} M \stackrel{\delta_-[-1]}{\longto} C_{\le 0} 
\]
of $M_-$ and
\[
M \stackrel{\pi}{\longto} M_+ \longto C_{\ge 1} \stackrel{\delta_+}{\longto} M[1]
\]
of $M_+$. 
\end{Rem}

Applying the theory to the examples from Section~\ref{1}, we get the following.

\begin{Def} \label{2C}
Assume that $k$ admits re\-so\-lution of singularities. Let $X \in Sm/k$.
Assume that the \emph{boundary motive} $\dMgm(X)$ of $X$ \cite[Def.~2.1]{W2}
admits a minimal weight filtration concentrated at $n=0$.
A Chow motive $M^{!*}(X)$, together with a factorization 
\[
\Mgm(X) \stackrel{i}{\longto} M^{!*}(X)
\stackrel{\pi}{\longto} \Mcgm(X)
\]
of the canonical morphism $u: \Mgm(X) \to \Mcgm(X)$ \cite[pp.~223--224]{V},
is called an \emph{absolute intersection motive of $X$} 
if it is a minimal factorization of $u$.
\end{Def}

Recall \cite[Prop.~2.2, Cor.~2.3]{W2} that $\dMgm(X)$
fits into a canonical exact triangle
\[
\dMgm(X) \longto \Mgm(X) \stackrel{u}{\longto} \Mcgm(X) \longto \dMgm(X)[1]
\]
in $\DgM_F \,$. 
Therefore, $C := \dMgm(X)[1]$ is a cone of $u$. 
According to \cite[Cor.~1.14]{W4},
the motive $\Mgm(X)$ belongs to $\DgM_{F,w \le 0} \,$, and the motive
with compact support $\Mcgm(X)$ to $\DgM_{F,w \ge 0} \,$.
Given Example~\ref{1Ba},
the assumptions of De\-fi\-nition~\ref{2B} are satisfied if
the boundary motive $\dMgm(X)$ belongs to the category $\DFAbgM \,$.

\begin{Def} \label{2D}
Assume $F$ to be a field contained in $\BR$, and let $X \in Sm/ \BC$.
An object $R \Gamma_{!*} (X,F(0))$ of $D^b ( \MHS_F )_{w = 0}$, together
with a factorization
\[
R \Gamma_c (X,F(0)) \stackrel{i}{\longto} R \Gamma_{!*} (X,F(0))
\stackrel{\pi}{\longto} R \Gamma (X,F(0))
\]
of the canonical morphism $u: R \Gamma_c (X,F(0)) \to R \Gamma (X,F(0))$,
is called an \emph{absolute Hodge theoretic intersection complex of $X$} 
if it is a minimal factorization of $u$.
\end{Def}

Here, the objects $R \Gamma_c (X,F(0))$ and $R \Gamma (X,F(0))$ of 
$D^b ( \MHS_F )$ are the canonical classes 
computing singular cohomology of $X(\BC)$
with and without compact support, respectively, 
and together with their Hodge structures \cite[Sect.~8.1]{D2}, \cite[Sect.~4]{Be}.
They satisfy the assumption on weights thanks to \cite[Cor.~(3.2.15)~(ii)]{D1}, \cite[Sect.~4.2]{Be}. 
According to Example~\ref{1I} and Corollary~\ref{1Ca}, any object of
$D^b ( \MHS_F )$ admits minimal weight filtrations.

\begin{Var} \label{2E}
(a)~Assume that $k$ admits re\-so\-lution of singularities. 
Let $X \in Sm/k$, and fix an idempotent endomorphism $e$ of the triangle
\[
\dMgm(X) \longto M(X) \longto \Mcgm(X) \longto \dMgm(X)[1] \; ,
\]
that is,
fix idempotent endomorphisms of each of the motives $M(X)$, $\Mcgm(X)$ and $\dMgm(X)$, 
which yield an endomorphism of the triangle.
Denote by $M(X)^e$, $\Mcgm(X)^e$ and $\dMgm(X)^e$ the images of $e$
on $M(X)$, $\Mcgm(X)$ and $\dMgm(X)$, respectively, 
and consider the canonical morphism
$u: M(X)^e \to \Mcgm(X)^e$. The object $M(X)^e$
belongs to $\DgM_{F,w \le 0} \,$, and 
$\Mcgm(X)^e$ to $\DgM_{F,w \ge 0} \, $. 
Assuming that $\dMgm(X)^e$ admits a minimal weight filtration concentrated at $n=0$,
Definition~\ref{2B} allows to define the notion of \emph{$e$-part of the
absolute intersection motive of $X$}, which is a triple $( M^{!*}(X)^e , i , \pi )$, with
$M^{!*}(X)^e \in \CHM_F$, and a minimal factorizaion
\[
\Mgm(X)^e \stackrel{i}{\longto} M^{!*}(X)^e
\stackrel{\pi}{\longto} \Mcgm(X)^e
\]
of $u$. The hypothesis on $\dMgm(X)^e$ is satisfied as soon as $\dMgm(X)^e \in \DFAbgM$ 
(which is the case in particular if the whole of $\dMgm(X)$ belongs to $\DFAbgM$). \\[0.1cm]
(b)~Similarly, for $F$ a field contained in $\BR$,
and $X \in Sm/\BC$, any pair $e$ of idempotent endomorphisms of 
$R \Gamma_c (X,F(0))$ and of $R \Gamma (X,F(0))$ commuting with the 
canonical morphism $R \Gamma_c (X,F(0)) \to R \Gamma (X,F(0))$
allows to define the notion of \emph{$e$-part of the
absolute intersection complex of $X$}, which is a 
triple $( R \Gamma_{!*} (X,F(0))^e , i , \pi )$, with
$R \Gamma_{!*} (X,F(0))^e \in D^b ( \MHS_F )_{w = 0}$, and a minimal factorizaion
\[
R \Gamma_c (X,F(0))^e \stackrel{i}{\longto} R \Gamma_{!*} (X,F(0))^e
\stackrel{\pi}{\longto} R \Gamma (X,F(0))^e
\]
of $u: R \Gamma_c (X,F(0))^e \to R \Gamma (X,F(0))^e$. \\[0.1cm]
(c)~In a similar vein, there are variants ``with coefficients'' of the absolute intersection
motive and the absolute Hodge theoretic intersection complex. 
The scheme $X$ need no longer be smooth,
and the constant coefficients $F(0)$ are replaced by a \emph{Chow motive} $N \in CHM(X)_F$
\cite[Def.~1.5]{W10} in the motivic context, and by a \emph{complex of Hodge modules}
$N \in D^b ( \MHM_F X )$ \cite[Sect.~4]{S}, which is pure of weight zero \cite[Sect.~4.5]{S}
in the Hodge theoretic context. Write $a$ for the structure morphism of $X$, and
consider the cano\-ni\-cal morphism $u: a_! (N) \to a_* (N)$ (\cite[Thm.~2.4.50~(2)]{CD}, 
\cite[proof of Thm.~4.3]{S}). 
The hypotheses on weights are satisfied (\cite[Thm.~3.8~($i_c$), ($i_c'$)]{He}, \cite[(4.5.2)]{S}). In the motivic
context, one needs to assume in addition that a cone of $u$ admits a minimal weight filtration
concentrated at $n=1$, which as before is the case if it belongs to $\DFAbgM \,$.
The result is an isomorphism class of a Chow motive $M^{!*}(X,N)$ in the motivic context,
and of an object $R \Gamma_{!*} (X,N)$ of $D^b ( \MHS_F )_{w = 0}$ in the Hodge theoretic
context, together with a minimal factorization of $u$. 
\end{Var}

Note that in the motivic context, following \cite[Cor.~16.1.6]{CD}, we have identified 
$\DgM_F$ with the triangulated category $\DBcFkM$ of \emph{constructible Beilinson motives}
over $\Spec k$ \cite[Def.~15.1.1]{CD}. Under that identification, the canonical morphism
$M(X) \to \Mcgm(X)$ equals the dual of the canonical morphism $a_! (\one_X) \to a_* (\one_X)$ 
(see \emph{e.g.} \cite[Prop.~2.8]{W12}). This point of view also allows to drop the
hypothesis on resolution of singularities in Definition~\ref{2C} and Variant~\ref{2E}~(a).

%%% Local Variables:
%%% mode: latex
%%% TeX-master: "head"
%%% End:

\bigskip
%\include{Sec3}
%%%%%%%%%%%%%%%%%%%%%%%%%%%%%%%%%%%%%%%%%%%%%%%%%%%%%%%%%%%%%%%%%%%%%%%
%
%  Section 3
%
%%%%%%%%%%%%%%%%%%%%%%%%%%%%%%%%%%%%%%%%%%%%%%%%%%%%%%%%%%%%%%%%%%%%%%%

\section{Absolute intersection cohomology}
\label{3}

%%%%%%%%%%%%%%%%%%%%%%%%%%%%%%%%

%%%%%%%%%%%%%%%%%%%%%%%%%%%%%%%%

The purpose of this section is to identify, up to isomorphism, the Hodge structure
on the cohomology objects of an absolute Hodge theoretic intersection complex 
$R \Gamma_{!*} (X,F(0))$ of a smooth, separated scheme $X$ over $\BC$. 
For $n \in \BZ$, denote by $H^n (X,F(0))$ the $n$-th singular cohomology group, and by
$H^n_c (X,F(0))$ the $n$-th cohomology group with compact support of $X(\BC)$. 

\begin{Def} \label{3A}
Assume $F$ to be a field contained in $\BR$, and let $X \in Sm/ \BC$.
\emph{Absolute intersection cohomology of $X$} is defined as the collection
\[
\bigl( H^n _{!*} (X,F(0)), i, \pi \bigr)_{n \in \BZ} \; ,
\]
where for $n \in \BZ$, we denote by $H^n _{!*} (X,F(0))$ the $n$-th cohomology
object of an absolute Hodge theoretic intersection complex 
$R \Gamma_{!*} (X,F(0))$ of $X$, and by 
\[
H^n_c (X,F(0)) \stackrel{i}{\longto} H^n_{!*} (X,F(0))
\stackrel{\pi}{\longto} H^n (X,F(0))
\]
the factorization of the canonical morphism $H^n_c (X,F(0)) \to H^n (X,F(0))$
induced by the factorization 
\[
R \Gamma_c (X,F(0)) \stackrel{i}{\longto} R \Gamma_{!*} (X,F(0))
\stackrel{\pi}{\longto} R \Gamma (X,F(0)) \; .
\]
\end{Def}

From the construction, we deduce the following.

\begin{Prop} \label{3B}
Assume $F$ to be a field contained in $\BR$, and let $X \in Sm/ \BC$. \\[0.1cm]
(a)~Absolute intersection cohomology of $X$ is equipped with a pure 
polarizable $F$-Hodge structure. More precisely, $H^n _{!*} (X,F(0))$ 
is pure (and polarizable) of weight $n$, for all $n \in \BZ$. \\[0.1cm]
(b)~Absolute intersection cohomology of $X$ is well-defined up to isomorphism
of $F$-Hodge structures. 
\end{Prop}

Remember that for $n \in \BZ$, the mixed graded-polarizable $F$-Hodge structures
on $H^n_c (X,F(0))$ and on $H^n (X,F(0))$ are of weights $\le n$ and $\ge n$, respectively.
Denote by $W_r$ the $r$-th filtration step
of the weight filtration of a mixed Hodge structure, and by $\Gr_r^W$ the quotient
of $W_r$ by $W_{r-1}$, $r \in \BZ$.

\begin{Cor} \label{3C}
Assume $F$ to be a field contained in $\BR$, let $X \in Sm/ \BC$ and $n \in \BZ$. 
Then $i: H^n_c (X,F(0)) \to H^n_{!*} (X,F(0))$ factors uniquely over 
\[
\Gr_n^W H^n_c (X,F(0)) = \frac{H^n_c (X,F(0))}{W_{n-1} H^n_c (X,F(0))} \; :
\]
\[
i: H^n_c (X,F(0)) \longonto \Gr_n^W H^n_c (X,F(0)) \stackrel{i_n}{\longto} H^n_{!*} (X,F(0)) \; ,
\]
and $\pi: H^n_{!*} (X,F(0)) \to H^n (X,F(0))$ factors uniquely over
\[
\Gr_n^W H^n (X,F(0)) = W_n H^n (X,F(0)) \; :
\]
\[
\pi: H^n_{!*} (X,F(0)) \stackrel{\pi_n}{\longto} \Gr_n^W H^n (X,F(0)) \longinto H^n (X,F(0)) \; .
\]
\end{Cor}

\begin{Def} \label{3D}
Assume $F$ to be a field contained in $\BR$, let $X \in Sm/ \BC$ and $n \in \BZ$. 
Denote by $u_n$ the canonical morphism 
\[
\Gr_n^W H^n_c (X,F(0)) \longto \Gr_n^W H^n (X,F(0))
\]
induced by $H^n_c (X,F(0)) \to H^n (X,F(0))$.
\end{Def}

Thus, any choice of absolute intersection cohomology of $X$ yields a factorization
\[
\Gr_n^W H^n_c (X,F(0)) \stackrel{i_n}{\longto} H^n_{!*} (X,F(0))
\stackrel{\pi_n}{\longto} \Gr_n^W H^n (X,F(0))
\]
of $u_n$, for all $n \in \BZ$. As we shall see, this factorization can be used to
charac\-terize absolute intersection cohomology up to isomorphism. \\

It turns out that the most appropriate context to formulate the result is purely
abstract. Fix a semi-simple Abelian category $\FA$.   

\begin{Def} \label{3E}
Let $v: S \to T$ be a morphism in $\FA$. Define
\[
\CH(v) := \ker(v) \oplus \imm(v) \oplus \coker(v) \; .
\]
\end{Def}

Monomorphisms and epimorphisms being split in $\FA$, 
\[
\CH(v) \cong S \oplus \coker(v) \quad \text{and} \quad \CH(v) \cong \ker(v) \oplus T \; ,
\]
the isomorphisms being in general non-canonical. In particular, $\CH(v) \cong S$
if $v$ is an epimorphism, and $\CH(v) \cong T$ if it is a monomorphism.
The inclusion of $\imm(v)$ into $\CH(v)$, or equivalently, the projection from $\CH(v)$
to $\imm(v)$ is an isomorphism if and only if $v$ is an isomorphism.

\begin{Prop} \label{3F}
Let $v: S \to T$ be a morphism in $\FA$. There exists a factorization
\[
S \stackrel{i^\CH}{\longto} \CH(v) \stackrel{\pi^\CH}{\longto} T
\]
of $v$ with the following properties: 
\begin{enumerate}
\item[(1)] $i^\CH$ is a monomorphism, and $\pi^\CH$ is an epimorphism.
\item[(2)] The restriction of $i^\CH$ to $\ker(v)$ equals the inclusion $i_1$
of the first component into $\CH(v) = \ker(v) \oplus \imm(v) \oplus \coker(v)$,
and the composition of $\pi^\CH$ with the quotient map to $\coker(v)$
equals the projection $\pi_3$ from $\CH(v) = \ker(v) \oplus \imm(v) \oplus \coker(v)$
to the last component. In particular,
\[
\pi_3 \circ i^\CH = 0: S \longto \coker(v) \quad \text{and} \quad 
\pi^\CH \circ i_1 = 0: \ker(v) \longto T \; .
\]
\end{enumerate}  
\end{Prop}

\begin{Proof}
Choose a left inverse $s$ of the inclusion of $\ker(v)$ into $S$, and a
right inverse $t$ of the quotient map from $T$ to $\coker(v)$. Define
\[
i^\CH := (s,q,0): S \longto \CH(v) = \ker(v) \oplus \imm(v) \oplus \coker(v) \; ,
\]
where $q$ is the quotient map from $S$ to $\imm(v)$, and
\[
\pi^\CH := 0 + \iota + t: \CH(v) = \ker(v) \oplus \imm(v) \oplus \coker(v) \longto T \; ,
\]
where $\iota$ is the inclusion of $\imm(v)$ into $T$. 
\end{Proof}

Again, factorizations as in Proposition~\ref{3F} are in general not unique
(they \emph{are} unique for trivial reasons, if $v$ is an isomorphism).
Still, any given choice satisfies a ``versal'' property, as we are about to see.

\begin{Thm} \label{3G}
Let $v: S \to T$ be a morphism in $\FA$. Fix a factorization
\[
S \stackrel{i^\CH}{\longto} \CH(v) \stackrel{\pi^\CH}{\longto} T
\]
of $v$ as in Proposition~\ref{3F}. Let
\[
S \stackrel{j}{\longto} H \stackrel{p}{\longto} T
\]
be any factorization of $v$, with a monomorphism $i$ and an epimorphism $p$.
Then there is an object $H'$ of $\FA$ and an isomorphism 
\[
H \cong \CH(v) \oplus H'
\]
compatible with the factorizations. In other words, denoting by $\iota$ the inclusion of,
and by $q$ the quotient map to $\CH(v)$, the diagrams
\[
\vcenter{\xymatrix@R-10pt{
        S \ar@{=}[d] \ar[r]^-{i^\CH} &  
        \CH(v) \ar@{^{ (}->}[d]_\iota \ar[r]^-{\pi^\CH} &
        T \ar@{=}[d] \\
        S \ar[r]^-j &  
        H \ar[r]^-p &
        T
\\}}
\]
and
\[
\vcenter{\xymatrix@R-10pt{
        S \ar@{=}[d] \ar[r]^-j &  
        H \ar@{>>}[d]_q \ar[r]^-p &
        T \ar@{=}[d] \\
        S \ar[r]^-{i^\CH} &  
        \CH(v) \ar[r]^-{\pi^\CH} &
        T
\\}}
\]
commute.
\end{Thm} 

\begin{Proof}
Denote by $s$ the composition
\[
S \stackrel{i^\CH}{\longto} \CH(v) \stackrel{\pi_1}{\longonto} \ker(v) \; ,
\]
where $\pi_1$ denotes the projection from $\CH(v)$ to the first component $\ker(v)$.
According to property~\ref{3F}~(2), the morphism $s$ is a left inverse of the 
inclusion of $\ker(v)$ into $S$. Let us show first that there is a left inverse of
the inclusion of $\ker(p)$ into $H$, which is compatible with $s$ in the sense that
the diagram
\[
\vcenter{\xymatrix@R-10pt{
        S \ar@{^{ (}->}[d]_j \ar[r]^-{s} &  
        \ker(v) \ar@{^{ (}->}[d]^j \\ 
        H \ar[r] &
        \ker(p)
\\}}
\]
commutes; we shall denote such a left inverse by the same letter $s$.

Choose a complement $H'$ of $j(\ker(v))$ in $\ker(p)$. Since $j^{-1}(\ker(p)) = \ker(v)$,
the intersection of $\imm(j)$ and $H'$ is trivial. Therefore,
\[
\ker(p) = j(\ker(v)) \oplus H' \longinto \imm(j) \oplus H' \subset H \; .
\]
Choose a complement $H''$ of $\imm(j) \oplus H'$ in $H$, and define
\[
s := (j(s),\id_{H'},0) : 
H = \imm(j) \oplus H' \oplus H'' \longto \ker(p) = j(\ker(v)) \oplus H' \oplus 0 \; .
\]
Then, using $s$, we have 
\[
\CH(v) = \ker(v) \oplus T \quad \text{and} \quad 
H = \ker(p) \oplus T \; ,
\]
$T$ being in both cases identified with $\ker(s)$. We define
\[
\iota:= (j,\id_T): \CH(v) \longto H 
\]
and 
\[
q:= (j^{-1} \circ \pi_1',\id_T): H \longto \CH(v) \; ,
\]
where $\pi_1'$ is the left inverse of $j: \ker(v) \into \ker(p)$ with
kernel equal to $H'$.
\end{Proof}

Let us get back to the geometric situation considered in the beginning of this section.

\begin{Thm} \label{3H}
Assume $F$ to be a field contained in $\BR$, let $X \in Sm/ \BC$ and $n \in \BZ$.
Recall the morphism $u_n: \Gr_n^W H^n_c (X,F(0)) \to \Gr_n^W H^n (X,F(0))$
(Definition~\ref{3D}). Fix a choice of 
$H^n_{!*} (X,F(0))$, together with the associated factorization
\[
\Gr_n^W H^n_c (X,F(0)) \stackrel{i_n}{\longto} H^n_{!*} (X,F(0))
\stackrel{\pi_n}{\longto} \Gr_n^W H^n (X,F(0))
\]
of $u_n$. \\[0.1cm]
(a)~The morphism $i_n$ is injective, and $\pi_n$ is surjective. \\[0.1cm]
(b)~In the semi-simple Abelian category of polarizable $F$-Hodge structures, 
which are pure of weight $n$, 
\[
H^n_{!*} (X,F(0)) \cong 
\CH(u_n) \; ,
\]
and the isomorphism can be chosen such that $i_n = i^\CH$ and $\pi_n = \pi^\CH$.
\end{Thm}

\begin{Proof}
According to Definition~\ref{3A}, absolute intersection cohomology is the
collection of cohomology objects of a choice of absolute Hodge theoretic
intersection complex $R \Gamma_{!*} (X,F(0))$. This choice induces a factorization
\[
R \Gamma_c (X,F(0)) \stackrel{i}{\longto} R \Gamma_{!*} (X,F(0))
\stackrel{\pi}{\longto} R \Gamma (X,F(0)) \; ,
\]
which on the level of cohomology yields the factorization
\[
H^n_c (X,F(0)) \longto H^n_{!*} (X,F(0)) \longto H^n (X,F(0)) \; .
\]
According to Definition~\ref{2D}, the above objects and morphisms fit into
a diagram
\[
\vcenter{\xymatrix@R-10pt{
        0 \ar[d] \ar[r] &
        C_{\ge 1}[-1] \ar[d] \ar@{=}[r] &
        C_{\ge 1}[-1] \ar[d]^{\delta[-1]} \ar[r] &
        0 \ar[d] \\
        R \Gamma_c (X,F(0)) \ar@{=}[d] \ar[r]^-i &  
        R \Gamma_{!*} (X,F(0)) \ar[d]^{\pi} \ar[r] &
        C_{\le 0} \ar[d] \ar[r] &
        R \Gamma_c (X,F(0))[1] \ar@{=}[d] \\
        R \Gamma_c (X,F(0)) \ar[d] \ar[r]^-u &
        R \Gamma (X,F(0)) \ar[d] \ar[r] &
        C \ar[d] \ar[r] &
        R \Gamma_c (X,F(0))[1] \ar[d] \\
        0 \ar[r] &
        C_{\ge 1} \ar@{=}[r] &
        C_{\ge 1} \ar[r] &
        0    
\\}}
\]
of exact triangles in $D^b ( \MHS_F )$,
where $\delta: C_{\ge 1} \to C_{\le 0}[1]$ belongs to the radical.

On the level of cohomology, the above induces a diagram of exact sequences,
the part of interest of which looks as follows.
\[
\vcenter{\xymatrix@R-10pt{
        H^{n-2}(C_{\ge 1}) \ar[d]^{H^{n-2}(\delta)} \ar[r] &
        0 \ar[d] \ar[r] &
        H^{n-1}(C_{\ge 1}) \ar[d] \ar@{=}[r] &
        H^{n-1}(C_{\ge 1}) \ar[d]^{H^{n-1}(\delta)} \\
        H^{n-1}(C_{\le 0}) \ar[d] \ar[r] &
        H^n_c (X,F(0)) \ar@{=}[d] \ar[r]^-i &  
        H^n_{!*} (X,F(0)) \ar[d]^{\pi} \ar[r] &
        H^n(C_{\le 0}) \ar[d] \\
        H^{n-1}(C) \ar[d] \ar[r] &
        H^n_c (X,F(0)) \ar[d] \ar[r] &
        H^n (X,F(0)) \ar[d] \ar[r] &
        H^n(C) \ar[d] \\
        H^{n-1}(C_{\ge 1}) \ar[r] &
        0 \ar[r] &
        H^n(C_{\ge 1}) \ar@{=}[r] &
        H^n(C_{\ge 1}) 
\\}}
\]
According to Corollary~\ref{1G}, the first and last columns equal the weight
filtrations of $H^{n-1}(C)$ and $H^n(C)$, concentrated at $n$ and $n+1$, respectively:
\[
\vcenter{\xymatrix@R-10pt{
        &
        &
        \frac{H^{n-1}(C)}{W_{n-1} H^{n-1}(C)} \ar[d] &
         \\
        W_{n-1} H^{n-1}(C) \ar[d] \ar[r] &
        H^n_c (X,F(0)) \ar@{=}[d] \ar[r]^-i &  
        H^n_{!*} (X,F(0)) \ar[d]^{\pi} \ar[r] &
        W_n H^n(C) \ar[d] \\
        H^{n-1}(C) \ar[r] &
        H^n_c (X,F(0)) \ar[r] &
        H^n (X,F(0)) \ar[d] \ar[r] &
        H^n(C) \ar[d] \\
        &
        &
        \frac{H^n(C)}{W_n H^n(C)} \ar@{=}[r] &
        \frac{H^n(C)}{W_n H^n(C)} 
\\}}
\]
Applying $\Gr^W_n$, the sequences remain exact, yielding
\[
\vcenter{\xymatrix@R-10pt{
        &
        &
        \Gr^W_n H^{n-1}(C) \ar[d] \\
        0 \ar[d] \ar[r] &
        \Gr^W_n H^n_c (X,F(0)) \ar@{=}[d] \ar[r]^-{i_n} &  
        H^n_{!*} (X,F(0)) \ar[d]^{\pi_n} \\
        \Gr^W_n H^{n-1}(C) \ar[r] &
        \Gr^W_n H^n_c (X,F(0)) \ar[r]^{u_n} &
        \Gr^W_n H^n (X,F(0)) \ar[d] \\
        &
        &
        0
\\}}
\] 
This diagram shows that $i_n$ is injective, and that $\pi_n$ is surjective,
proving part~(a) of our claim.

The last column also shows that $\ker(\pi_n)$ equals the image of $\Gr^W_n H^{n-1}(C)$
in $H^n_{!*} (X,F(0))$. The morphism $\Gr^W_n H^{n-1}(C) \to H^n_{!*} (X,F(0))$
factors over $\Gr^W_n H^n_c (X,F(0))$. The image of 
$\Gr^W_n H^{n-1}(C) \to \Gr^W_n H^n_c (X,F(0))$ equals $\ker(u_n)$,
according to the second line. The morphism $i_n$ being injective,
we get a canonical short exact sequence
\[
0 \longto \ker(u_n) \longto H^n_{!*} (X,F(0)) \stackrel{\pi_n}{\longto} \Gr^W_n H^n (X,F(0))
\longto 0 \; ,
\] 
which is (in general non-canonically) split. This shows that
\[
H^n_{!*} (X,F(0)) \cong \ker(u_n) \oplus \Gr^W_n H^n (X,F(0)) \; ,
\]
which in turn is isomorphic to $\CH(u_n)$. Applying Theorem~\ref{3G} to
\[
\Gr_n^W H^n_c (X,F(0)) \stackrel{i_n}{\longto} H^n_{!*} (X,F(0))
\stackrel{\pi_n}{\longto} \Gr_n^W H^n (X,F(0))
\]
(using part~(a)), we see that the isomorphism $H^n_{!*} (X,F(0)) \isoto \CH(u_n)$
can be chosen in a way compatible with the factorizations.
\end{Proof}

\begin{Cor} \label{3I}
Assume $F$ to be a field contained in $\BR$, let $X \in Sm/ \BC$ and $n \in \BZ$.
Then \emph{interior cohomology} $H^n_! (X,F(0))$ is a direct factor of any choice
of absolute intersection cohomology $H^n_{!*} (X,F(0))$. The two are cano\-ni\-cally
isomorphic if and only if $u_n: \Gr^W_n H^n_c (X,F(0)) \to \Gr^W_n H^n (X,F(0))$
is an isomorphism.
\end{Cor} 

\begin{Proof}
By definition, $H^n_! (X,F(0)) = \imm(u_n)$,
and $\CH(u_n)$ contains $\imm(u_n)$ as a direct factor, with complement equal to
$\ker(u_n) \oplus \coker(u_n)$.
Now apply Theorem~\ref{3H}~(b). 
\end{Proof}

Note that a sufficient condition for $u_n$ to be an isomorphism is that
\emph{boun\-dary cohomology} $\partial H^r (X,F(0))$ avoids weight $n$ in degrees
$r = n-1$ and $r = n$.

\begin{Rem} 
It happens rarely that $\partial H^r (X,F(0))$ avoids weights $r$ and $r+1$
in all degrees $r$. Note however that all constructions and results concerning
Hodge structures on cohomology established in this section admit obvious analogues
in the context considered in Variant~\ref{2E}~(b), \emph{i.e.}, in the presence
of an idempotent endomorphism $e$, and that
there exist non-trivial examples where  
$\partial H^r (X,F(0))^e$ \emph{does} avoid weights $r$ and $r+1$
in all degrees.
\end{Rem}

\begin{Cor} \label{3K}
Assume $F$ to be a field contained in $\BR$, let $X \in Sm/ \BC$ and $\widetilde{X}$
any smooth compactification of $X$.
Then absolute intersection cohomology $H^n_{!*} (X,F(0))$ is a direct factor
of $H^n (\widetilde{X},F(0))$, for all $n \in \BZ$.
\end{Cor} 

\begin{Proof}
The morphism $u_n$ factors canonically through $H^n (\widetilde{X},F(0))$:
\[
\Gr^W_n H^n_c (X,F(0)) \stackrel{j}{\longto} H^n (\widetilde{X},F(0)) 
\stackrel{p}{\longto} \Gr^W_n H^n (X,F(0)) \; .
\]
According to \cite[Cor.~(3.2.17)]{D1}, the morphism $p$ is surjective. Therefore, by duality,
the morphism $j$ is injective.
Now apply Theorems~\ref{3H}~(b) and \ref{3G}.
\end{Proof}

\begin{Rem} \label{3L}
Both complexes $R \Gamma_{!*} (X,F(0))$ and $R \Gamma (\widetilde{X},F(0))$ being 
pure of weight zero, they are (in general, non-canonically) isomorphic to the direct
sum of their $n$-th cohomology objects, shifted by $n$. According to Corollary~\ref{3K},
$R \Gamma_{!*} (X,F(0))$ is therefore isomorphic to
a direct factor of $R \Gamma (\widetilde{X},F(0))$.
We do not know whether the isomorphism can be chosen in a way compatible with the
factorizations of 
\[
R \Gamma_c (X,F(0)) \longto R \Gamma (X,F(0)) 
\]
through $R \Gamma_{!*} (X,F(0))$ and through $R \Gamma (\widetilde{X},F(0))$
(see Remark~\ref{2Ba}~(a)).
\end{Rem}

%%% Local Variables:
%%% mode: latex
%%% TeX-master: "head"
%%% End:

\bigskip
%\include{Sec4}
%%%%%%%%%%%%%%%%%%%%%%%%%%%%%%%%%%%%%%%%%%%%%%%%%%%%%%%%%%%%%%%%%%%%%%%
%
%  Section 4
%
%%%%%%%%%%%%%%%%%%%%%%%%%%%%%%%%%%%%%%%%%%%%%%%%%%%%%%%%%%%%%%%%%%%%%%%

\section{Examples}
\label{4}

%%%%%%%%%%%%%%%%%%%%%%%%%%%%%%%%

%%%%%%%%%%%%%%%%%%%%%%%%%%%%%%%%

Let us begin by recalling a special case of a result of Voevodsky's.

\begin{Thm} \label{4A}
Let $r$ and $s$ be two integers. Then
\[
\Hom_{\DgM_F} \bigl( \BZ(r)[2r], \BZ(s)[2s] \bigr) = 0 \quad \text{if} \quad r \ne s \; .
\]
\end{Thm} 

\begin{Proof}
This is a special case of \cite[Cor.~2]{V2} (\cite[Prop.~4.2.9]{V} if $k$
admits resolution of singula\-ri\-ties).
\end{Proof}

In the sequel, we assume that the base field $k$ admits resolution of singularities.
For $X \in Sm/k$, repeat Construction~\ref{Cons} for the canonical morphism
$u: \Mgm(X) \to \Mcgm(X)$ and (the shift by $+1$ of) a weight filtration 
\[
\dm_{\le -1} \longto \dMgm(X) \longto \dm_{\ge 0} 
\stackrel{d}{\longto} \dm_{\le -1}[1]
\]
of $\dMgm(X)$
(with $\dm_{\le -1} \in \CC_{w \le -1}$ and $\dm_{\ge 0} \in \CC_{w \ge 0}$):
\[
\vcenter{\xymatrix@R-10pt{
        0 \ar[d] \ar[r] &
        \dm_{\ge 0} \ar[d]_{d_+} \ar@{=}[r] &
        \dm_{\ge 0} \ar[d]^{d} \ar[r] &
        0 \ar[d] \\
        \Mgm(X) \ar@{=}[d] \ar[r]^-i &  
        M \ar[d]_{\pi} \ar[r]^-{d_-} &
        \dm_{\le -1}[1] \ar[d] \ar[r] &
        \Mgm(X)[1] \ar@{=}[d] \\
        \Mgm(X) \ar[d] \ar[r]^-u &
        \Mcgm(X) \ar[d] \ar[r] &
        \dMgm(X)[1] \ar[d] \ar[r] &
        \Mgm(X)[1] \ar[d] \\
        0 \ar[r] &
        \dm_{\ge 0}[1] \ar@{=}[r] &
        \dm_{\ge 0}[1] \ar[r] &
        0    
\\}}
\]
The motive $M$ is pure of weight zero, \emph{i.e.}, it is a Chow motive, 
\[
\Mgm(X) \stackrel{i}{\longto} M \stackrel{\pi}{\longto} \Mcgm(X)
\]
is a factorization of $u$, and 
\[
\dm_{\le -1} \longto \Mgm(X) \stackrel{i}{\longto} M 
\stackrel{d-}{\longto} \dm_{\le -1}[1]
\]
and
\[
M \stackrel{\pi}{\longto} \Mcgm(X) \longto \dm_{\ge 0}[1]
\stackrel{d_+[1]}{\longto} M[1]
\]
are weight filtrations of $\Mgm(X)$ and of $\Mcgm(X)$, respectively.
Proposition~\ref{2A} tells us that up to isomorphism,
the process is reversible: a triple $(M,i,\pi)$ leads to a weight filtration
$(\dm_{\le -1},\dm_{\ge 0})$ of $\dMgm(X)$. \\

According to Definition~\ref{2C}, the triple $(M,i,\pi)$ is 
an absolute intersection motive of $X$ (in which case we write $M = M^{!*}(X)$) if 
and only if it is a minimal factorization of $u$, \emph{i.e.}, if and only if the morphism 
\[
\dm_{\ge 0} 
\stackrel{d}{\longto} \dm_{\le -1}[1]
\]
belongs to the radical $\rad_{\DgM_F} ( \dm_{\ge 0} , \dm_{\le -1}[1] )$.

\begin{Ex} \label {4C}
Put $X = \BA^1_k$. The morphism $u: \Mgm(\BA^1_k) \to \Mcgm(\BA^1_k)$ factors canonically
through the Chow motive $\Mgm(\BP^1_k)$:
\[
u: \Mgm(\BA^1_k) \stackrel{j_*}{\longto} \Mgm(\BP^1_k) = \Mcgm(\BP^1_k) 
\stackrel{j^*}{\longto} \Mcgm(\BA^1_k) \; ,
\]
where $j^*$ and $j_*$ denote the morphisms induced by the open immersion $j$ of
$\BA^1_k$ into $\BP^1_k$. We claim that $(\Mgm(\BP^1_k),j_*,j^*)$ is a minimal
factorization of $u$.

Indeed, both $j^*$ and $j_*$ can be completed to give weight filtrations
of $\Mgm(\BA^1_k)$ and of $\Mcgm(\BA^1_k)$, respectively: the inclusion of
the point at infinity of $\Mgm(\BP^1_k)$ yields exact (\emph{purity}
and \emph{localization}) triangles
\[
\BZ(1)[1] \longto \Mgm(\BA^1_k) \stackrel{j_*}{\longto} \Mgm(\BP^1_k) 
\stackrel{d_-}{\longto} \BZ(1)[2] 
\]
\cite[Prop.~3.5.4]{V} and
\[
\Mgm(\BP^1_k) \stackrel{j^*}{\longto} \Mcgm(\BA^1_k) \longto \BZ(0)[1]
\stackrel{d_+[1]}{\longto} \Mgm(\BP^1_k)[1]
\]
\cite[Prop.~4.1.5]{V}. We thus get a diagram of exact triangles
\[
\vcenter{\xymatrix@R-10pt{
        0 \ar[d] \ar[r] &
        \BZ(0) \ar[d]_{d_+} \ar@{=}[r] &
        \BZ(0) \ar[d]^{d} \ar[r] &
        0 \ar[d] \\
        \Mgm(\BA^1_k) \ar@{=}[d] \ar[r]^-{j_*} &  
        \Mgm(\BP^1_k) \ar[d]_{j^*} \ar[r]^-{d_-} &
        \BZ(1)[2] \ar[d] \ar[r] &
        \Mgm(\BA^1_k)[1] \ar@{=}[d] \\
        \Mgm(\BA^1_k) \ar[d] \ar[r]^-u &
        \Mcgm(\BA^1_k) \ar[d] \ar[r] &
        \dMgm(\BA^1_k)[1] \ar[d] \ar[r] &
        \Mgm(\BA^1_k)[1] \ar[d] \\
        0 \ar[r] &
        \BZ(0)[1] \ar@{=}[r] &
        \BZ(0)[1] \ar[r] &
        0    
\\}}
\]
The morphism $d: \BZ(0) \to \BZ(1)[2]$ is zero according to Theorem~\ref{4A}.
\emph{A fortiori}, it belongs to 
$\rad_{\DgM_F} ( \BZ(0), \BZ(1)[2] )$.
Therefore, $(\Mgm(\BP^1_k),j_*,j^*)$ is minimal; in particular, 
\[
M^{!*}(\BA^1_k) = \Mgm(\BP^1_k) \; .
\]
The third column of the above diagram also shows that
\[
\dMgm(\BA^1_k) \cong \BZ(1)[1] \oplus Z(0) \; .
\]
\end{Ex}

\begin{Exo} \label {4D}
More generally, if $X$ is the complement of a zero-dimensio\-nal sub-scheme $Z$
in a smooth, proper $k$-scheme $Y$, which is of pure dimension $n \ge 1$, then  
\[
M^{!*}(X) = \Mgm(Y) 
\]
and
\[
\dMgm(X) \cong \Mgm(Z)(n)[2n-1] \oplus \Mgm(Z) \; .
\]
\end{Exo}

\begin{Ex} \label{4E}
To generalize further, assume that $X$ is the complement of a smooth closed 
sub-scheme $Z$ in a smooth, proper $k$-scheme $Y$, and that the immersion
of $Z$ into $Y$ is of pure codimension $c$.
Again, we have weight filtrations
of $\Mgm(X)$ and of $\Mcgm(X)$, respectively, coming from
purity and localization:
\[
\Mgm(Z)(c)[2c - 1] \longto \Mgm(X) \stackrel{j_*}{\longto} \Mgm(Y) 
\stackrel{d_-}{\longto} \Mgm(Z)(c)[2c]
\]
\cite[Prop.~3.5.4]{V} and
\[
\Mgm(Y) \stackrel{j^*}{\longto} \Mcgm(X) \longto \Mgm(Z)[1]
\stackrel{d_+[1]}{\longto} \Mgm(Y)[1]
\]
\cite[Prop.~4.1.5]{V}. We get
\[
\vcenter{\xymatrix@R-10pt{
        0 \ar[d] \ar[r] &
        \Mgm(Z) \ar[d]_{d_+} \ar@{=}[r] &
        \Mgm(Z) \ar[d]^{d} \ar[r] &
        0 \ar[d] \\
        \Mgm(X) \ar@{=}[d] \ar[r]^-{j_*} &  
        \Mgm(Y) \ar[d]_{j^*} \ar[r]^-{d_-} &
        \Mgm(Z)(c)[2c] \ar[d] \ar[r] &
        \Mgm(X)[1] \ar@{=}[d] \\
        \Mgm(X) \ar[d] \ar[r]^-u &
        \Mcgm(X) \ar[d] \ar[r] &
        \dMgm(X)[1] \ar[d] \ar[r] &
        \Mgm(X)[1] \ar[d] \\
        0 \ar[r] &
        \Mgm(Z)[1] \ar@{=}[r] &
        \Mgm(Z)[1] \ar[r] &
        0    
\\}}
\]
The (shift by $-1$ of the) third column is clearly a weight filtration
of the boundary motive $\dMgm(X)$; 
the question is to determine whether it is minimal, \emph{i.e.},
whether the morphism
\[
\Mgm(Z) \stackrel{d}{\longto} \Mgm(Z)(c)[2c]
\]
belongs to the radical $\rad_{\DgM_F} ( \Mgm(Z) , \Mgm(Z)(c)[2c] )$. 
Observe that according to \cite[Thm.~4.3.7~3., Prop.~4.2.9]{V},
\[
\Hom_{\DgM_F} \bigl( \Mgm(Z) , \Mgm(Z)(c)[2c] \bigr) = CH^{dim_X}(Z \times_k Z)_F \; .
\]
(a)~Assume that the dimension of $Z$ is strictly smaller than $c$. Then
\[
CH^{dim_X}(Z \times_k Z)_F = 0 \; .
\]
Therefore, the morphism $d$ is zero, $(\Mgm(Y),j_*,j^*)$ is minimal, 
\[
M^{!*}(X) = \Mgm(Y) \; ,
\]
and
\[
\dMgm(X) \cong \Mgm(Z)(c)[2c-1] \oplus \Mgm(Z) \; .
\]
These four statements remain true, up to replacing $\Mgm(Z)(c)[2c-1]$
by $\Mgm(Z)^*(dim_X)[2dim_X-1]$, if the smoothness assumption on $Z$ is dropped. 
(Hint: use induction on the dimension of $T \in Sch/k$ to show that for all non-negative integers $j$,
\[
\Hom_{\DgM_F} ( \Mgm(T) , \BZ(r)[2r-j] ) = 0 \quad \text{if} \quad r > dim_T + j \; .
\]
For the induction step, use \cite[Prop.~4.1.3]{V}.
Apply this to $T = Z \times_k Z$, $j=0$ and $r = dim_X$.) \\[0.1cm]
(b)~By contrast, without the condition on the codimension from (a), the morphism
\[
\Mgm(Z) \stackrel{d}{\longto} \Mgm(Z)(c)[2c] \; ,
\]
which corresponds to a class in the Chow group
\[
CH^{dim_X}(Z \times_k Z)_F \; ,
\]
can \emph{a priori} be non-zero. Actually, it does happen that $d$ does not belong to the radical,
meaning that $(\Mgm(Y),j_*,j^*)$ is not minimal. 
This is the case for $X = \BA^n_k$ and $Y = \BP^n_k$, $n \ge 2$; 
we refer to Example~\ref{4G} below. \\[0.1cm]
(c)~Given the definition of composition of cycles, we see that
under the identification
\[
\Hom_{\DgM_F} \bigl( \Mgm(Z) , \Mgm(Z)(c)[2c] \bigr) = CH^{dim_X}(Z \times_k Z)_F \; ,
\]
the morphism $d$ corresponds to the image of the class 
\[
[\Delta_Z] \in CH^{dim_X}(Z \times_k Y)_F
\]
of the diagonal $\Delta_Z$ under the pull-back 
\[
CH^{dim_X}(Z \times_k Y)_F \longto CH^{dim_X}(Z \times_k Z)_F \; . 
\]
Equivalently (see \emph{e.g.} \cite[Lemma~1.1~(i)]{Ku}), it corresponds to the 
\emph{Lefschetz operator} $L_{i^*[Z]}$ associated to the ``self intersection''
$i^*[Z]$, defined as the image under 
\[
\Delta_{Z,*} : CH^c(Z)_F \longto CH^{dim_X}(Z \times_k Z)_F
\]
of $i^*[Z]$, where $[Z]$ denotes the class of the cycle $Z$ in $CH^c(Y)_F$,
and $i^*$ the pull-back from $CH^c(Y)_F$ to $CH^c(Z)_F \,$. \\[0.1cm]
(d)~Assume that $c=1$, \emph{i.e.}, that $Z$ is a smooth divisor. The morphism
\[
\Mgm(Z) \stackrel{d}{\longto} \Mgm(Z)(1)[2] 
\]
corresponds to 
\[
L_{i^*[Z]} \in CH^{dim_X}(Z \times_k Z)_F \; .
\]
Assume that $i^*[Z]$ or its opposite is ample. 
Assume also $\Mgm(Z)$ satisfies the following \emph{weak form of the Lefschetz decomposition}
for $L_{i^*[Z]}$: there exist two decompositions
\[
\Mgm(Z) = P_1 \oplus \Mgm(Z)^s_1 \quad \text{and} \quad \Mgm(Z) = P_2 \oplus \Mgm(Z)^s_2 \; , 
\]  
such that $L_{i^*[Z]}$ equals zero on $P_1$, maps the sub-motive $\Mgm(Z)^s_1$ of $\Mgm(Z)$ to 
the sub-motive $\Mgm(Z)^s_2(1)[2]$ of $\Mgm(Z)(1)[2]$, and
induces an isomorphism $\Mgm(Z)^s_1 \isoto \Mgm(Z)^s_2(1)[2]$.

This means that the zero morphism $P_1 \to P_2(1)[2]$ and $d$ have isomorphic cones. Thus,
\[
P_2(1)[1] \longto \dMgm(X) \longto P_1 \stackrel{0}{\longto} P_2(1)[2]
\]
is a minimal weight filtration of $\dMgm(X)$. The absolute intersection motive 
$M^{!*}(X)$ can be constructed out of $\Mgm(Y)$ as follows. Start with the diagram
\[
\vcenter{\xymatrix@R-10pt{
        0 \ar[d] \ar[r] &
        \Mgm(Z) \ar[d]_{d_+} \ar@{=}[r] &
        \Mgm(Z) \ar[d]^{d} \ar[r] &
        0 \ar[d] \\
        \Mgm(X) \ar@{=}[d] \ar[r]^-{j_*} &  
        \Mgm(Y) \ar[d]_{j^*} \ar[r]^-{d_-} &
        \Mgm(Z)(1)[2] \ar[d] \ar[r] &
        \Mgm(X)[1] \ar@{=}[d] \\
        \Mgm(X) \ar[d] \ar[r]^-u &
        \Mcgm(X) \ar[d] \ar[r] &
        \dMgm(X)[1] \ar[d] \ar[r] &
        \Mgm(X)[1] \ar[d] \\
        0 \ar[r] &
        \Mgm(Z)[1] \ar@{=}[r] &
        \Mgm(Z)[1] \ar[r] &
        0    
\\}}
\]
The motives $\Mgm(Z)^s_1$ and $\Mgm(Z)^s_2(1)[2]$ are
direct factors of $\Mgm(Z)$ and $\Mgm(Z)(1)[2]$, respectively, and the morphism
$d = d_- \circ d_+$ induces an isomorphism between them. Therefore, the kernel $M_{gm}'(Y)$ of
the composition
\[
\Mgm(Y) \stackrel{d_-}{\longto} \Mgm(Z)(1)[2] \longonto \Mgm(Z)^s_2(1)[2]
\]
(exists and) is a direct factor of $\Mgm(Y)$, with complement equal to the image of
the composition
\[
\Mgm(Z)^s_1 \longinto \Mgm(Z) \stackrel{d_+}{\longto} \Mgm(Y) \; .
\]
By definition, the morphisms $j_*$ and $j^*$ factor through $M_{gm}'(Y)$.
We thus get the following direct factor of the above diagram.
\[
\vcenter{\xymatrix@R-10pt{
        0 \ar[d] \ar[r] &
        P_1 \ar[d] \ar@{=}[r] &
        P_1 \ar[d]^{0} \ar[r] &
        0 \ar[d] \\
        \Mgm(X) \ar@{=}[d] \ar[r]^-{j_*} &  
        M_{gm}'(Y) \ar[d]_{j^*} \ar[r] &
        P_2(1)[2] \ar[d] \ar[r] &
        \Mgm(X)[1] \ar@{=}[d] \\
        \Mgm(X) \ar[d] \ar[r]^-u &
        \Mcgm(X) \ar[d] \ar[r] &
        \dMgm(X)[1] \ar[d] \ar[r] &
        \Mgm(X)[1] \ar[d] \\
        0 \ar[r] &
        P_1[1] \ar@{=}[r] &
        P_1[1] \ar[r] &
        0    
\\}}
\] 
The zero morphism belonging to the radical, we see that 
$(M_{gm}'(Y),j_*,j^*)$ is an absolute intersection motive of $X$. In particular, 
\[
M^{!*}(X) = M_{gm}'(Y) \; .
\]
\end{Ex}

\begin{Rem}
(a)~Recall the classical notion of \emph{Lefschetz decomposition for $L_D$}
associated to an ample divisor $D$ on a smooth, proper $k$-scheme $Z$
of pure dimension $r$: 
there is a decomposition
\[
\Mgm(Z) = \bigoplus_{i=0}^{2r} \bigoplus_{m = max(0,i-r)}^{[\frac{i}{2}]} L^m P^{i-2m}
\]
into Chow motives $L^m P^{i-2m}$, such that \\[0.1cm]
(1)~for all $0 \le m \le r-i-1$, the morphism $L_D$ induces an isomorphism
between the sub-motive $L^m P^i$ of $\Mgm(Z)$ and
the sub-motive $L^{m+1} P^i(1)[2]$ of $\Mgm(Z)(1)[2]$, \\[0.1cm]
(2)~for all $i \le r$, the morphism $L_D$ is zero on $L^{r-i} P^i$. \\[0.1cm]
Such decompositions exist if $Z$ is an Abelian variety \cite[Thm.~5.1]{Ku}.  

Putting 
\[
P_1 := \bigoplus_{i=0}^{r} L^{r-i} P^i \quad \text{and} \quad 
P_2 := \bigoplus_{i=0}^{r} P^i
\]
(and $\Mgm(Z)^s_1$, $\Mgm(Z)^s_2$ equal to the sums of the respective remaining
direct factors in the Lefschetz decomposition), we get 
\[
\Mgm(Z) = P_1 \oplus \Mgm(Z)^s_1 \quad \text{and} \quad \Mgm(Z) = P_2 \oplus \Mgm(Z)^s_2 \; , 
\]  
such that $L_D$ equals zero on $P_1$, and
induces an isomorphism $\Mgm(Z)^s_1 \isoto \Mgm(Z)^s_2(1)[2]$. \\[0.1cm]
(b)~We leave it to the reader to formulate a generalization of Example~\ref{4E}~(d)
for smooth sub-schemes $i: Z \into Y$ of pure codimension $c$, such that
$i^*[Z]$ is the $c$-fold self intersection of the class of an ample divisor on $Z$. 
\end{Rem}

Let us come back to the situation considered in Example~\ref{4E}.

\begin{Thm} \label{4F}
Let $Z$ be a smooth, geometrically connected closed sub-scheme of pure dimension $c \ge 1$ 
of a smooth, proper $k$-scheme $Y$ of pure dimension $2c$.
Denote by $i$ the closed immersion of $Z$,
and by $j$ the open immersion of its complement $X$ into $Y$.
Consider the morphism 
\[
\Mgm(Z) \stackrel{d}{\longto} \Mgm(Z)(c)[2c] 
\]
corresponding to 
\[
L_{i^*[Z]} \in CH^{2c}(Z \times_k Z)_F = CH_0 (Z \times_k Z)_F \; .
\]
(a)~The morphism $d$ belongs to 
$\rad_{\DgM_F} ( \Mgm(Z) , \Mgm(Z)(c)[2c] )$ if and only if the self
intersection number $Z \cdot Z$ of $Z$ in $Y$ equals zero. \\[0.1cm]
(b)~If $Z \cdot Z \ne 0$, then there exist unique decompositions
\[
\Mgm(Z) = M_1^r \oplus M_1^s \; , \; \Mgm(Z)(c)[2c] = M_2^r \oplus M_2^s \; , 
\]  
such that 
\begin{enumerate}
\item[(1)] the decompositions are respected by $d$:
\[
d = d^r \oplus d^s \in \Hom_{\DgM_F} \bigl( \Mgm(Z) , \Mgm(Z)(c)[2c] \bigr) \; ,
\]
\item[(2)] the morphism $d^r$ belongs to $\rad_{\DgM_F} (M_1^r,M_2^r)$, 
\item[(3)] and the morphism $d^s$ is an isomorphism $M_1^s \isoto M_2^s$.
\end{enumerate} 
Furthermore, the motives $M_1^s$ and $M_2^s$ are canonically isomorphic to
$\BZ(c)[2c]$, and under these isomorphisms, 
\begin{enumerate}
\item[(4)] the composition $\BZ(c)[2c] \isoto M_1^s \into \Mgm(Z)$ corresponds to
the inclusion into $\Mgm(Z)$ of the canonical sub-motive $M_{gm}^{2c}(Z)$ \cite[Sect.~1.13]{Sch},
\item[(5)] the composition $\Mgm(Z)(c)[2c] \onto M_2^s \isoto \BZ(c)[2c]$ corresponds to
the twist by $c$ of the shift by $2c$ of
the projection from $\Mgm(Z)$ onto the canonical quotient $M_{gm}^{0}(Z)$ \cite[Sect.~1.11]{Sch},
\item[(6)] and the morphism $d^s$ corresponds to multiplication by $Z \cdot Z$.
\end{enumerate} 
\end{Thm}

\begin{Proof}
We have
\[
\Hom_{\DgM_F} \bigl( \Mgm(Z)(c)[2c] , \Mgm(Z) \bigr) = CH^0(Z \times_k Z)_F 
\]
\cite[Thm.~4.3.7~3., Prop.~4.2.9]{V}. 
The group $CH^0(Z \times_k Z)$ is generated by the class of $Z \times_k Z$.
According to \cite[Lemma~1.1~(i)]{Ku}, the composition $e'$ of the latter with
$L_{i^*[Z]}$ equals $Z \cdot Z$ times an idempotent endomorphism of $\Mgm(Z)$.
In particular, $\id_{\Mgm(Z)} - e'$ is an automorphism if and only if 
$e' = 0$, \emph{i.e.}, if and only if $Z \cdot Z = 0$.

This shows part~(a) of our claim. As for part~(b), assume that $Z \cdot Z \ne 0$.
Put $f := \frac{1}{Z \cdot Z} [Z \times_k Z] \in CH^0(Z \times_k Z)_F$. 
As the reader will verify, both
\[
e:= f \circ d \in \End_{\DgM_F} \bigl( \Mgm(Z) \bigr)
\]
and 
\[
g:= d \circ f \in \End_{\DgM_F} \bigl( \Mgm(Z)(c)[2c] \bigr)
\]
are idempotent; more precisely, the image $M_1^s$ of $e$ equals $M_{gm}^{2c}(Z)$, and
the image $M_2^s$ of $f$ projects isomorphically onto $M_{gm}^{0}(Z)(c)[2c]$. 
Put 
\[
M_1^r := \ker(e) \quad \text{and} \quad M_2^r := \ker(f) \; .
\] 
Properties~(1)--(6) then follow from our construction, and from the fact that
any morphism in $\Hom_{\DgM_F} ( \Mgm(Z)(c)[2c] , \Mgm(Z))$ is a scalar multiple
of $f$ --- hence 
\[
\Hom_{\DgM_F} \bigl( M_2^r , M_1^r \bigr) = 0 \; .
\]
The same argument shows the unicity of the decompositions.
\end{Proof}

\begin{Cor} \label{4Fa}
In the situation of Theorem~\ref{4F}, assume that $Z \cdot Z = 0$.
Then $(\Mgm(Y),j_*,j^*)$ is an absolute intersection motive of $X$. In particular, 
\[
M^{!*}(X) = \Mgm(Y) \; .
\]
\end{Cor}

\begin{Rem}
The above setting provides examples for non-zero morphisms in the radical.
Namely, if $i^*[Z]$ is of degree zero, without being rationally equivalent to zero,
then $0 \ne d : \Mgm(Z) \to \Mgm(Z)(c)[2c]$ belongs to the radical. 
In this case, 
\[
\dMgm(X) \not\cong \Mgm(Z)(c)[2c-1] \oplus \Mgm(Z) \; .
\]
\end{Rem}

\begin{Cor} \label{4Fb}
In the situation of Theorem~\ref{4F}, assume that $Z \cdot Z \ne 0$. \\[0.1cm]
(a)~There is a canonical direct factor $M_{gm}'(Y)$ of $\Mgm(Y)$, admitting a 
cano\-ni\-cal complement, which is isomorphic to $\BZ(c)[2c]$. \\[0.1cm]
(b)~The morphism $j_*: \Mgm(X) \to \Mgm(Y)$ factors through the sub-motive $M_{gm}'(Y)$,
and the morphism $j^*: \Mgm(Y) \to \Mcgm(X)$ factors through the quotient $M_{gm}'(Y)$. \\[0.1cm]
(c)~The triplet
$(M_{gm}'(Y),j_*,j^*)$ is an absolute intersection motive of $X$. In particular, 
\[
\Mgm(Y) \cong \BZ(c)[2c] \oplus M^{!*}(X) \; .
\]
\end{Cor}

\begin{Proof}
Recall the diagram
\[
\vcenter{\xymatrix@R-10pt{
        0 \ar[d] \ar[r] &
        \Mgm(Z) \ar[d]_{d_+} \ar@{=}[r] &
        \Mgm(Z) \ar[d]^{d} \ar[r] &
        0 \ar[d] \\
        \Mgm(X) \ar@{=}[d] \ar[r]^-{j_*} &  
        \Mgm(Y) \ar[d]_{j^*} \ar[r]^-{d_-} &
        \Mgm(Z)(c)[2c] \ar[d] \ar[r] &
        \Mgm(X)[1] \ar@{=}[d] \\
        \Mgm(X) \ar[d] \ar[r]^-u &
        \Mcgm(X) \ar[d] \ar[r] &
        \dMgm(X)[1] \ar[d] \ar[r] &
        \Mgm(X)[1] \ar[d] \\
        0 \ar[r] &
        \Mgm(Z)[1] \ar@{=}[r] &
        \Mgm(Z)[1] \ar[r] &
        0    
\\}}
\]
According to Theorem~\ref{4F}~(b), the motive $\BZ(c)[2c]$ can be identified with
direct factors of both $\Mgm(Z)$ and $\Mgm(Z)(c)[2c]$, and the morphism
$d = d_- \circ d_+$ induces an isomorphism between them. Therefore, the kernel $M_{gm}'(Y)$ of
the composition
\[
\Mgm(Y) \stackrel{d_-}{\longto} \Mgm(Z)(c)[2c] \longonto \BZ(c)[2c]
\]
(exists and) is a direct factor of $\Mgm(Y)$, with complement equal to the image of
the composition
\[
\BZ(c)[2c] \longinto \Mgm(Z) \stackrel{d_+}{\longto} \Mgm(Y) \; .
\]
By definition, the morphisms $j_*$ and $j^*$ factor through $M_{gm}'(Y)$.
With the notations of Theorem~\ref{4F}~(b), we thus get the following direct factor of the above diagram.
\[
\vcenter{\xymatrix@R-10pt{
        0 \ar[d] \ar[r] &
        M_1^r \ar[d] \ar@{=}[r] &
        M_1^r \ar[d]^{d^r} \ar[r] &
        0 \ar[d] \\
        \Mgm(X) \ar@{=}[d] \ar[r]^-{j_*} &  
        M_{gm}'(Y) \ar[d]_{j^*} \ar[r] &
        M_2^r \ar[d] \ar[r] &
        \Mgm(X)[1] \ar@{=}[d] \\
        \Mgm(X) \ar[d] \ar[r]^-u &
        \Mcgm(X) \ar[d] \ar[r] &
        \dMgm(X)[1] \ar[d] \ar[r] &
        \Mgm(X)[1] \ar[d] \\
        0 \ar[r] &
        M_1^r[1] \ar@{=}[r] &
        M_1^r[1] \ar[r] &
        0    
\\}}
\] 
But according to Theorem~\ref{4F}~(b), 
\[
d^r \in \rad_{\DgM_F} \left( M_1^r,M_2^r \right) \; .
\]
\end{Proof}

\begin{Ex} \label{4G}
Let $n \ge 1$, and put $X = \BA^n_k$. We have $\Mgm(\BA^n_k) = \BZ(0)$
and $\Mcgm(\BA^n_k) = \BZ(n)[2n]$ \cite[Cor.~4.1.8]{V}.
In particular, the exact triangle
\[
\Mcgm(\BA^n_k)[-1] \longto \dMgm(\BA^n_k) \longto \Mgm(\BA^n_k)
\stackrel{u}{\longto} \Mcgm(\BA^n_k)
\]
is a weight filtration of $\dMgm(\BA^n_k)$.
By Theorem~\ref{4A}, the morphism $u$
is zero; therefore, the above weight filtration is minimal. Construction~\ref{Cons}
ensures the existence of a Chow motive $M$, and of morphisms $i$, $\pi$, $d_+$
and $d_-$ such that 
\[
\vcenter{\xymatrix@R-10pt{
        0 \ar[d] \ar[r] &
        \Mgm(\BA^n_k) \ar[d]_{d_+} \ar@{=}[r] &
        \Mgm(\BA^n_k) \ar[d]^{0} \ar[r] &
        0 \ar[d] \\
        \Mgm(\BA^n_k) \ar@{=}[d] \ar[r]^-i &  
        M \ar[d]_{\pi} \ar[r]^-{d_-} &
        \Mcgm(\BA^n_k) \ar[d] \ar[r] &
        \Mgm(\BA^n_k)[1] \ar@{=}[d] \\
        \Mgm(\BA^n_k) \ar[d] \ar[r]^-0 &
        \Mcgm(\BA^n_k) \ar[d] \ar[r] &
        \dMgm(\BA^n_k)[1] \ar[d] \ar[r] &
        \Mgm(\BA^n_k)[1] \ar[d] \\
        0 \ar[r] &
        \Mgm(\BA^n_k)[1] \ar@{=}[r] &
        \Mgm(\BA^n_k)[1] \ar[r] &
        0    
\\}}
\]
is a diagram of exact sequences.
Observe that both morphisms $\Mcgm(\BA^n_k) \to \Mgm(\BA^n_k)[1]$ in the second line 
and in the second row are zero
since their source is of weight zero, and their target, of weight one.
Thus, we see that $M:= \Mgm(\BA^n_k) \oplus \Mcgm(\BA^n_k)$, $i:= d_+ :=$ the inclusion
of the first direct factor, and $\pi := d_- :=$ the projection to the second direct
factor, is a solution. Therefore, $(\Mgm(\BA^n_k) \oplus \Mcgm(\BA^n_k),i,\pi)$ is minimal; in particular, 
\[
M^{!*}(\BA^n_k) = \Mgm(\BA^n_k) \oplus \Mcgm(\BA^n_k) = \BZ(0) \oplus \BZ(n)[2n] \; ;
\]
we may think of this as the ``motive of the $2n$-sphere'' $\Mgm(\BS^{2n})$. 
The minimal weight filtration of $\dMgm(\BA^n_k)$ shows that
\[
\dMgm(\BA^n_k) \cong \Mcgm(\BA^n_k)[-1] \oplus \Mgm(\BA^n_k) = \BZ(n)[2n-1] \oplus \BZ(0) \; .
\]
\end{Ex}

\begin{Rem} \label{4Ga}
Example~\ref{4G} shows in particular that in general, the absolute intersection motive
of a product is unequal to the tensor product of the absolute intersection motives
of the factors.
\end{Rem}

\begin{Ex} \label{4H}
Let $X^*$ be a proper surface over $k$, which we assume to be normal.
The singular locus $X^*_{sing}$ of $X^*$ is of dimension zero; denote by $X$ its complement.
We claim that the absolute intersection motive of $X$ (exists and) is isomorphic
to the \emph{intersection motive of $X^*$} \cite{CM,W7}.

Choose a resolution of singularities. More precisely, consider in addition the following diagram,
assumed to be cartesian:
\[
\vcenter{\xymatrix@R-10pt{
        X \ar@{^{ (}->}[r]^-{j} \ar@{=}[d] &
        Y \ar@{<-^{ )}}[r]^-{i} \ar[d]_\pi &
        Z \ar[d]^\pi \\
        X \ar@{^{ (}->}[r] &
        X^* \ar@{<-^{ )}}[r] &
        X^*_{sing}
\\}}
\]
where $\pi$ is proper and birational, $Y$ is smooth (and proper), and $Z$ is a divisor with
normal crossings, whose irreducible components $Z_m$ are smooth 
(this is possible according to \cite[Theorem]{L2} and the discussion in \cite[pp.~191--194]{L1}).

As before, localization gives a weight filtration of $\Mcgm(X)$:
\[
\Mgm(Y) \stackrel{j^*}{\longto} \Mcgm(X) \longto \Mgm(Z)[1]
\stackrel{d_+[1]}{\longto} \Mgm(Y)[1] \; .
\]
Dualizing it, twisting by $2$ and shifting by $4$, we get a weight filtration of $\Mgm(X)$,
according to \cite[Thm.~4.3.7~3.]{V}:
\[
\Mgm(Z)^*(2)[3] \longto \Mgm(X) \stackrel{j_*}{\longto} \Mgm(Y) 
\stackrel{d_-}{\longto} \Mgm(Z)^*(2)[4] \; .
\]
We get
\[
\vcenter{\xymatrix@R-10pt{
        0 \ar[d] \ar[r] &
        \Mgm(Z) \ar[d]_{d_+} \ar@{=}[r] &
        \Mgm(Z) \ar[d]^{d} \ar[r] &
        0 \ar[d] \\
        \Mgm(X) \ar@{=}[d] \ar[r]^-{j_*} &  
        \Mgm(Y) \ar[d]_{j^*} \ar[r]^-{d_-} &
        \Mgm(Z)^*(2)[4] \ar[d] \ar[r] &
        \Mgm(X)[1] \ar@{=}[d] \\
        \Mgm(X) \ar[d] \ar[r]^-u &
        \Mcgm(X) \ar[d] \ar[r] &
        \dMgm(X)[1] \ar[d] \ar[r] &
        \Mgm(X)[1] \ar[d] \\
        0 \ar[r] &
        \Mgm(Z)[1] \ar@{=}[r] &
        \Mgm(Z)[1] \ar[r] &
        0    
\\}}
\]
There is a canonical morphism
\[
\iota_*: \bigoplus_m M^2_{gm}(Z_m) \longinto \bigoplus_m \Mgm(Z_m) \longto \Mgm(Z)
\]
\cite[Sect.~1.13]{Sch}. According to \cite[Thm.~2.2~(i)]{W7}
(see also \cite[Sect.~2.5]{CM}), the composition
\[
\bigoplus_m M^2_{gm}(Z_m) \stackrel{\iota_*}{\longto} \Mgm(Z) 
\stackrel{d}{\longto} \Mgm(Z)^*(2)[4] 
\stackrel{(\iota_*)^*(2)[4]}{\longto} \bigoplus_m M^2_{gm}(Z_m)^*(2)[4]
\] 
is an isomorphism. Therefore, the kernel of
the composition
\[
\Mgm(Y) \stackrel{d_-}{\longto} \Mgm(Z)^*(2)[4] 
\stackrel{(\iota_*)^*(2)[4]}{\longto} \bigoplus_m M^2_{gm}(Z_m)^*(2)[4]
\]
(exists and) is a direct factor of $\Mgm(Y)$, with complement equal to the image of
the composition
\[
\bigoplus_m M^2_{gm}(Z_m) \stackrel{\iota_*}{\longto} \Mgm(Z) 
\stackrel{d_+}{\longto} \Mgm(Y) \; .
\]
But by definition \cite[Def.~2.3, Ex.~5.2]{W7}, that kernel equals $M^{!*}(X^*)$,
the intersection motive of $M^*$. 
The morphisms $j_*$ and $j^*$ factor through $M^{!*}(X^*)$.
Denoting by $M_{gm}^{\le 1}(Z)$ the complement of $\oplus_m M^2_{gm}(Z_m)$ in $\Mgm(Z)$
(see \cite[Lemma~5.4]{W7}),
we thus get the following direct factor of the above diagram.
\[
\vcenter{\xymatrix@R-10pt{
        0 \ar[d] \ar[r] &
        M_{gm}^{\le 1}(Z) \ar[d] \ar@{=}[r] &
        M_{gm}^{\le 1}(Z) \ar[d]^{d^{\le 1}} \ar[r] &
        0 \ar[d] \\
        \Mgm(X) \ar@{=}[d] \ar[r]^-{j_*} &  
        M^{!*}(X^*) \ar[d]_{j^*} \ar[r] &
        M_{gm}^{\le 1}(Z)^*(2)[4] \ar[d] \ar[r] &
        \Mgm(X)[1] \ar@{=}[d] \\
        \Mgm(X) \ar[d] \ar[r]^-u &
        \Mcgm(X) \ar[d] \ar[r] &
        \dMgm(X)[1] \ar[d] \ar[r] &
        \Mgm(X)[1] \ar[d] \\
        0 \ar[r] &
        M_{gm}^{\le 1}(Z)[1] \ar@{=}[r] &
        M_{gm}^{\le 1}(Z)[1] \ar[r] &
        0    
\\}}
\]
Here, we set $d^{\le 1}:=$ the restriction of $d$ to $M_{gm}^{\le 1}(Z)$.
We claim that
\[
d^{\le 1} \in \rad_{\DgM_F} \left( M_{gm}^{\le 1}(Z) , M_{gm}^{\le 1}(Z)^*(2)[4] \right) \; .
\]
In fact, \emph{all} morphisms $M_{gm}^{\le 1}(Z) \to M_{gm}^{\le 1}(Z)^*(2)[4]$ belong
to the radical, since
\[
\Hom_{\DgM_F} \left( M_{gm}^{\le 1}(Z)^*(2)[4] , M_{gm}^{\le 1}(Z)  \right) = 0 \; .
\]
To prove this latter claim, observe that the closed covering by the $Z_m$
induces a weight filtration of $M_{gm}^{\le 1}(Z)$ (see the first part of the proof 
of \cite[Prop.~6.5~(i)]{W7}):
\[
\bigoplus_{n<m} \Mgm(Z_{n,m}) \longto \bigoplus_m M_{gm}^{\le 1}(Z_m) \longto M_{gm}^{\le 1}(Z)
\longto \bigoplus_{n<m} \Mgm(Z_{n,m})[1] \; ,
\]
where we denote by $Z_{n,m}$ the intersection $Z_n \cap Z_m$, and by $M_{gm}^{\le 1}(Z_m)$ the 
quotient of $\Mgm(Z_m)$ by
$M^2_{gm}(Z_m)$. Whence the dual weight filtration
\[
\bigoplus_{n<m} \Mgm(Z_{n,m})^*[-1] \longto M_{gm}^{\le 1}(Z)^* \longto
\bigoplus_m M_{gm}^{\le 1}(Z_m)^* \longto \bigoplus_{n<m} \Mgm(Z_{n,m})^*
\]
of $M_{gm}^{\le 1}(Z)^*$. Orthogonality for weight structures shows that
any morphism from $M_{gm}^{\le 1}(Z)^*(2)[4]$ to a motive of non-negative weights
factors through $\oplus_m M_{gm}^{\le 1}(Z_m)^*(2)[4]$, and that any morphism from
a motive of weight zero to $M_{gm}^{\le 1}(Z)$ factors through $\oplus_m M_{gm}^{\le 1}(Z_m)$.
Therefore,
\[
\Hom_{\DgM_F} \left( M_{gm}^{\le 1}(Z)^*(2)[4] , M_{gm}^{\le 1}(Z)  \right) 
\]
is a direct factor of
\[
\bigoplus_{n,m} \Hom_{\DgM_F} \left( M_{gm}^{\le 1}(Z_n)^*(2)[4] , M_{gm}^{\le 1}(Z_m)  \right) \; . 
\]
We claim that
\[
\Hom_{\DgM_F} \left( M_{gm}^{\le 1}(Z_n)^*(2)[4] , M_{gm}^{\le 1}(Z_m)  \right)  = 0  
\]
for all $n$ and $m$. Indeed, the split monomorphisms $M^2_{gm}(Z_m) \into \Mgm(Z_m)$
induce an isomorphism between 
\[
\Hom_{\DgM_F} \left( M_{gm}^2(Z_n)^*(2)[4] , M_{gm}^2(Z_m)  \right) 
\]
and
\[
\Hom_{\DgM_F} \left( \Mgm(Z_n)^*(2)[4] , \Mgm(Z_m)  \right) \; ,
\]
as can be seen from the comparison wih Chow theory: both groups are identified with 
\[
CH^0(Z_n \times_k Z_m)_F 
\]
\cite[Thm.~4.3.7~3., Prop.~4.2.9]{V}, and the morphism between them corresponds
to the identity. 

Altogether, we proved that $(M^{!*}(X^*),j_*,j^*)$ is minimal.
In particular,
\[
M^{!*}(X) = M^{!*}(X^*) \; .
\]
\end{Ex}

For $k = \BC$, and $F$ a field contained in $\BR$, the reader may choose to
compute absolute intersection cohomology in the geometric situations treated 
in this section,
taking into account the results from Section~\ref{3}, in particular, Theorem~\ref{3H}.
The computations are \emph{a priori} compatible with the above under the Hodge theoretic
realization.

%%% Local Variables:
%%% mode: latex
%%% TeX-master: "head"
%%% End:

\bigskip
%\include{Sec5}
%%%%%%%%%%%%%%%%%%%%%%%%%%%%%%%%%%%%%%%%%%%%%%%%%%%%%%%%%%%%%%%%%%%%%%%
%
%  Section 5
%
%%%%%%%%%%%%%%%%%%%%%%%%%%%%%%%%%%%%%%%%%%%%%%%%%%%%%%%%%%%%%%%%%%%%%%%

\section{A question}
\label{5}

%%%%%%%%%%%%%%%%%%%%%%%%%%%%%%%%

%%%%%%%%%%%%%%%%%%%%%%%%%%%%%%%%

In this section, the base $k$ equals the field $\BC$ of complex numbers.
The coefficients $F$ are assumed to be contained in $\BR$.
Let $X \in Sm/\BC$.
The examples treated in Section~\ref{4}
seem to suggest the following.

\begin{Ques} \label{5A}
Does there exist an open, dense immersion $j$ of $X(\BC)$ into a compact topological space $X^{!*}$
satisfying the hypothesis of \cite[2.1.16]{BBD}, and equipped with a stratification into topological manifolds
of even (real) dimension, among which $X(\BC)$, such that 
``intersection cohomology of $X^{!*}$ equals absolute intersection cohomology of $X$''
in the following sense:
for all $n \in \BZ$, there is an isomorphism of $F$-vector spaces
\[
H^n (X^{!*},\ujast F) \isoto H^n_{!*} (X,F)
\]
from cohomology of $X^{!*}$ with coefficients in the intermediate extension $\ujast F$
to the $F$-vector space $H^n_{!*} (X,F)$ underlying 
absolute intersection cohomology $H^n_{!*} (X,F(0))$, respecting the factorizations 
\[
H^n_c (X,F) \stackrel{H^n j_*}{\longto} H^n (X^{!*},\ujast F)
\stackrel{H^n j^*}{\longto} H^n (X,F)
\]
and 
\[
H^n_c (X,F(0)) \stackrel{i}{\longto} H^n_{!*} (X,F(0))
\stackrel{\pi}{\longto} H^n (X,F(0)) \; ?
\]
\end{Ques}

A word of explanation is in order. Since $X$ is smooth, the constant sheaf $F$ on $X$ 
can be seen as a perverse sheaf up to a shift by $-dim_X$:
\[
F = (F[dim_X])[-dim_X] \; ,
\]
where $F[dim_X]$ belongs to the heart of the perverse $t$-structure. 
By slight abuse of notation, we define
\[
\ujast F := (\ujast F[dim_X])[-dim_X] \; .
\]
Note that $\ujast F[dim_X]$ is defined, and 
that if $X^{!*}$ is a topological manifold, then $H^n (X^{!*},\ujast F)$ equals
$H^n (X^{!*},F)$ for all $n \in \BZ$ \cite[Prop.~2.1.17]{BBD}. Note also that
\[
H^n j_*: H^n_c (X,F) \longto H^n (X^{!*},\ujast F)
\]
factors through the $F$-vector space underlying $\Gr_n^W H^n_c (X,F(0))$, and
\[
H^n j^*: H^n (X^{!*},\ujast F) \longto H^n (X,F)
\]
factors through the $F$-vector space underlying $\Gr_n^W H^n (X,F(0))$ 
(Corol\-la\-ry~\ref{3C}).

\begin{Exs} \label{5B}
The answer to Question~\ref{5A} is affirmative in the following cases. \\[0.1cm]
(a)~For $X = \BA^1_\BC$, put $X^{!*} = \BP^1(\BC)$ (Example~\ref{4C}). \\[0.1cm]
(b)~More generally,
for $X$ equal to a finite number of closed points
in a smooth, proper $\BC$-scheme $Y$, which is of pure dimension $n \ge 1$,
put $X^{!*} = Y(\BC)$ (Example~\ref{4D}). \\[0.1cm]
(c)~Even more generally,
let $X$ be equal to the complement of a closed 
sub-scheme $Z$ in a smooth, proper $\BC$-scheme $Y$. Assuming that the dimension of $Z$
is strictly smaller than its codimension in $Y$, we may put
$X^{!*} = Y(\BC)$ (Example~\ref{4E}~(a)). \\[0.1cm]
(d)~Let $X$ be equal to the complement of a smooth closed 
sub-scheme $Z$ of pure dimension $c \ge 1$
in a smooth, proper $\BC$-scheme $Y$ of pure dimension $2c$.
Assume that $Z$ is connected, and that the self intersection number
$Z \cdot Z$ of $Z$ in $Y$ is non-zero. 
We know (Corollary~\ref{4Fb}) that $X^{!*}$ cannot be chosen to be equal
to $Y(\BC)$. \\[0.1cm] 
(e)~Let $X$ be equal to the complement of a smooth divisor
$Z$ in a smooth, proper $\BC$-scheme $Y$.
Assume that the pull-back $i^*[Z]$ of the divisor to the sub-scheme $Z$ is ample, 
and that $\Mgm(Z)$ satisfies the weak form of the Lefschetz decomposition
for $L_{i^*[Z]}$. 
We know (Example~\ref{4E}~(d)) that $X^{!*}$ cannot be chosen to be equal
to $Y(\BC)$. \\[0.1cm] 
(f)~For $X = \BA^n_\BC$, $n \ge 1$,
put $X^{!*} = \BS^{2n}$ (Example~\ref{4G}). \\[0.1cm]
(g)~For $X$ equal to the complement of the singular locus in
a proper normal surface $X^*$ over $\BC$, we may put
$X^{!*} = X^*(\BC)$ (Example~\ref{4H}).
\end{Exs}

We omitted to specify $X^{!*}$ in Examples~\ref{5B}~(d) and (e). Before doing so
(see Proposition~\ref{5D}), let us treat the remaining example
from Section~\ref{4}; we claim that it is fundamentally different from 
those treated in Example~\ref{5B}.

\begin{Ex} \label{5C}
Let $X$ be equal to the complement of a smooth closed 
sub-scheme $Z$ of pure dimension $c \ge 1$
in a smooth, proper $\BC$-scheme $Y$ of pure dimension $2c$.
Assume that $Z$ is connected, and that 
$Z \cdot Z = 0$. We may then put
$X^{!*} = Y(\BC)$ (Corollary~\ref{4Fa}). 
\end{Ex}

Denote by $X^+$ the Alexandrov one-point compactification of $X(\BC)$.

\begin{Prop} \label{5D}
(a)~In all of Examples~\ref{5B}, the compactification $X^{!*}$
of $X(\BC)$ may be chosen to be equal to $X^+$. \\[0.1cm]
(b)~In Example~\ref{5C}, exactly one of the following cases occurs: either
\[
dim_F H^n (X^+,\ujast F) > dim_F H^n_{!*} (X,F) \quad \text{for} \quad 
n = 2c-1 \; \text{and} \; n = 2c+1  
\]
or there is no isomorphism
\[ 
H^{2c} (X^+,\ujast F) \longto H^{2c}_{!*} (X,F)
\]
respecting the factorizations of 
$H^{2c}_c (X,F) \to H^{2c} (X,F)$.
In particular, the compactification $X^{!*}$
of $X(\BC)$ may not be chosen to be equal to $X^+$.
\end{Prop}

Proposition~\ref{5D} is obvious in certain cases. For the treatment of the others,
the following observation turns out to be useful.

\begin{Prop} \label{5Da}
Let $X \in Sm/\BC$ of pure dimension. Denote by $j$ the immersion of
$X(\BC)$ into $X^+$. Then
\[
H^n (X^+,\ujast F) = H^n (X, F) \quad \text{for} \quad n \le  dim_X - 1 \; ,
\]
\[
H^{dim_X} (X^+,\ujast F) = H^{dim_X}_! (X, F) \; ,
\]
and
\[
H^n (X^+,\ujast F) = H^n_c (X, F) \quad \text{for} \quad n \ge dim_X + 1 \; .
\]
Each of the identifications is understood with the canonical factorization 
of the morphism $H^n_c (X, F) \to H^n (X, F)$.
\end{Prop}

\begin{Rem} \label{5Db}
(a)~In particular, intersection cohomology of $X^+$ carries a canonical Hodge structure
(which in general is mixed). \\[0.1cm]
(b)~Given the formulae from Proposition~\ref{5Da}, we see that in those cases where
intersection cohomology of $X^+$ equals absolute intersection cohomology of $X$,
the isomorphisms $H^n (X^+,\ujast F) \isoto H^n_{!*} (X,F)$, $n \in \BZ$ are unique,
and respect the Hodge structures (which are pure).
\end{Rem}

\begin{Proofof}{Proposition~\ref{5Da}}
Denote by $\tau^{t \le \bullet}$ 
and $\tau^{t \ge \bullet}$ 
the truncation functors with respect to the $t$-structure
on the derived category of sheaves on the complement $\{\ast\}$ of $X(\BC)$ in $X^+$,
in other words, on the derived category of $F$-vector spaces.
Denote by $i$ the inclusion of $\{\ast\}$ into $X^+$.
According to \cite[Prop.~2.1.11]{BBD}, there is an exact triangle
\[
i_* \tau^{t \ge dim_X} i^* R j_* F [-1] \longto \ujast F \longto R j_* F \longto
i_* \tau^{t \ge dim_X} i^* R j_* F \; .
\]
Applying the direct image $a_*$ under the structure map $a$ of $X^+$, we get
isomorphisms
\[
H^n (X^+,\ujast F) \isoto H^n (X, F) \quad \text{for} \quad n < dim_X \; ,
\]
and a monomorphism
\[
H^{dim_X} (X^+,\ujast F) \longinto H^{dim_X} (X, F) \; .
\]
Dually, we obtain an epimorphism
\[
H^{dim_X}_c (X, F)  \longonto H^{dim_X} (X^+,\ujast F) \; ,
\]
and isomorphisms
\[
H^n_c (X, F) \isoto H^n (X^+,\ujast F) \quad \text{for} \quad n > dim_X \; .
\]
\end{Proofof}

Using Proposition~\ref{5Da}, the reader may choose
to provide a direct proof of Proposition~\ref{5D}. The one we shall give is a consequence
of the next result.

\begin{Thm} \label{5E}
Let $X \in Sm/\BC$ of pure dimension. Then the following are equivalent.
\begin{enumerate}
\item[(1)] The compactification $X^{!*}$
of $X(\BC)$ (exists and) may be chosen to be equal to $X^+$.
\item[(2)] For all integers $n \le dim_X-1$,
the Hodge structure on $\partial H^n (X,F(0))$ is of weights $\le n$.
\item[(3)] For all integers $n \ge dim_X$,
the Hodge structure on $\partial H^n (X,F(0))$ is of weights $\ge n+1$.  
\item[(4)] There exists an open immersion of $X$ into a smooth, proper $\BC$-scheme $Y$,
with complement $i: Z \into Y$, such that the maps
\[
H^n(Z, i^! F) \longto H^n(Z,F)
\]
($i^! = $ the exceptional inverse image under the immersion $i$)
are injective in the range $2 \cdot codim_Y Z \le n \le dim_X$.  
\item[(5)] Whenever $X$ is represented 
as the complement of a closed 
sub-scheme $i: Z \into Y$ in a smooth, proper $\BC$-scheme $Y$, the maps
\[
H^n(Z, i^! F) \longto H^n(Z,F)
\]
are injective in the range $2 \cdot codim_Y Z \le n \le dim_X$.   
\end{enumerate} 
If the equivalent conditions~(1)--(5) are satisfied, then the following also hold.
\begin{enumerate}
\item[(6)] For all integers $n \le dim_X-1$,
the Hodge structure on $H^n (X,F(0))$ is pure of weight $n$.
\item[(7)] For all integers $n \ge dim_X+1$,
the Hodge structure on $H^n_c (X,F(0))$ is pure of weight $n$.  
\end{enumerate} 
\end{Thm} 

\begin{Proof}
Proposition~\ref{5Da} shows that claim~(1) implies claims~(6) and (7). 

The Hodge structures on $H^n_c (X,F(0))$ and on
$H^{2dim_X - n} (X,F(0))(dim_X)$ are dual to each other. Given the long exact sequence
\[
\ldots \partial H^{n-1} (X,F(0)) \longto H^n_c (X,F(0)) \stackrel{u}{\longto} H^n (X,F(0))
\longto \partial H^n (X,F(0)) \; ,
\] 
claims~(2) and (3) are therefore equivalent. 

We leave it as an exercice to the reader to show, using Theorem~\ref{3H}~(b)
and Proposition~\ref{5Da}, that (2) and (3) together (hence individually)
are equivalent to (1).
 
To finish the proof, note that we have nothing to prove if $X$
is proper. Else, let $i: Z \into Y$ be a closed immersion, with $Y$ smooth and
proper over $\BC$, and such that $X = Y - Z$. Write $j$ for the open immersion
of $X$ into $Y$, and $a: Y \to \Spec \BC$ for the structure morphism of $Y$. 
The morphism $u: R \Gamma_c (X,F(0)) \to R \Gamma (X,F(0))$ in $D^b ( \MHS_F )$
is the result of applying the direct image $a_*$ to the morphism
\[
v: j_! F(0) \longto j_* F(0)
\]
in $D^b ( \MHM_F Y )$ \cite[proof of Thm.~4.3]{S}. A canonical choice of cone of $v$ is given
by $i_* i^* j_* F(0)$ (apply \cite[(4.4.1)]{S} to $j_* F(0)$). Thus, $a_* i_* i^* j_* F(0)$ is a cone of $u$;
its cohomology objects are thus equal to boundary cohomology.

Consider the localization exact triangle
\[
i_* i^! F(0) \longto F(0) \longto j_* F(0) \longto i_* i^! F(0)[1]
\]
in $D^b ( \MHM_F Y )$ \cite[(4.4.1)]{S}. Applying the functor $a_* i_* i^*$ yields
\[
a_* i_* i^! F(0) \longto a_* i_* i^* F(0) \longto a_* i_* i^* j_* F(0) \longto 
a_* i_* i^! F(0)[1] \; .
\]
We thus see that boundary cohomology is part of a long exact sequence
\[
\ldots \partial H^{n-1} (X,F(0)) \longto H^n (Z,i^!F(0)) \longto H^n (Z,F(0))
\longto \partial H^n (X,F(0)) \; .
\] 
The Hodge structures on $H^n (Z,i^!F(0))$ and on $H^n (Z,F(0))$
are of weights $\ge n$ and $\le n$, respectively \cite[(4.5.2)]{S}.
Therefore, for a fixed integer $n$,
the Hodge structure on $\partial H^n (X,F(0))$ is of weights $\le n$
if and only if $H^{n+1} (Z,i^!F) \to H^{n+1} (Z,F)$ is injective.
We leave it to the reader to show, for example by using duality, that
\[
H^n (Z,i^!F) = 0 \quad \text{for} \quad n < 2 \cdot codim_Y Z \; .
\]
\end{Proof}

\begin{Ex}
For the surface $X = \BA^1_\BC \times_\BC \BP^1_\BC$, 
the Hodge structure on $H^n (X,F(0))$ is pure of weight $n$ for $n = 0, 1$,
meaning that property~\ref{5E}~(6) is satisfied.
But according to Proposition~\ref{5D}, the compactification $X^{!*}$
of $X(\BC)$ may not be chosen to be equal to $X^+$.
\end{Ex}

\begin{Proofof}{Proposition~\ref{5D}}
The claim is obvious in cases~\ref{5B}~(a) and \ref{5B}~(f). Cases~\ref{5B}~(b) and \ref{5B}~(g)
follow from invariance of intersection cohomology
under direct images of finite morphisms.

To treat cases~\ref{5B}~(c), \ref{5B}~(d) and \ref{5C}, we apply criteria
\ref{5E}~(4) (for \ref{5B}~(c) and \ref{5B}~(d)) and \ref{5E}~(5) (for \ref{5C})
to the respective choice of $i:Z \into Y$. In case~\ref{5B}~(c), it is trivially satisfied
since $2 \cdot codim_Y Z > dim_X$. In cases~\ref{5B}~(d) and \ref{5C}, the only
degree that needs to be verified is $n = 2c$. By purity, the 
source $H^{2c}(Z, i^! F_Y)$ of the map in question is identified with 
$H^0(Z,F(-c)) = F(-c)$. The trace map gives the same type of identification for the target
$H^{2c}(Z,F)$. Under these identifications,  
\[
H^{2c}(Z, i^! F_Y) \longto H^{2c}(Z,F)
\]
corresponds to multiplication by $Z \cdot Z$. 
According to Theorem~\ref{5E}, 
we may put $X^{!*} = X^+$ in case~\ref{5B}~(d), but not in case~\ref{5C}.
In the latter case, note that by duality, 
\[
dim_F H^{2c-1} (X^+,\ujast F) > dim_F H^{2c-1}_{!*} (X,F) 
\]
if and only if
\[
dim_F H^{2c+1} (X^+,\ujast F) > dim_F H^{2c+1}_{!*} (X,F) \; , 
\]
and recall that for all integers $n$, 
\[
H^n_{!*} (X,F) = H^n (Y,F) \; .
\]
The sequence
\[
H^{2c}(Z,F) \stackrel{\partial}{\longto} H^{2c+1}_c(X,F) = H^{2c+1} (X^+,\ujast F)
\longto H^{2c+1}(Y,F) \longto 0
\]
is exact. Thus, if $\partial \ne 0$, then 
$dim_F H^{2c+1} (X^+,\ujast F) > dim_F H^{2c+1}_{!*} (X,F)$.
If $\partial = 0$, then 
\[
H^{2c}_c(X,F) \longto H^{2c}(Y,F)
\]
is not surjective. This means that interior cohomology 
\[
H^{2c}_!(X,F) = H^{2c} (X^+,\ujast F)
\]
and $H^{2c}(Y,F)$ cannot be identified in a way compatible with 
the factorizations of $H^{2c}_c (X, F) \to H^{2c} (X, F)$.

To treat case~\ref{5B}~(e), note that the maps
\[
H^n(Z, i^! F(0)) \longto H^n(Z,F(0)) 
\]
are identified with the cohomological Lefschetz operators
\[
L_{i^*[Z]} : H^{n-2}(Z, F(-1)) \longto H^n(Z,F(0)) \; .
\] 
According to the weak Lefschetz theorem, these are injective for all $n \le dim_Z+1 = dim_X$.
In other words, criterion~(4) from Theorem~\ref{5E} is satisfied.  
\end{Proofof}

\begin{Exo} \label{5F}
For $X = \BG_{m,\BC} \times \BA^1_\BC$, compute $M^{!*}(X)$ and identify a choice
of $X^{!*}$. Prove that $X^{!*}$ can neither be chosen to be equal to $X^+$,
nor to the space of $\BC$-valued points of a smooth algebraic compactification of $X$. 
Relate this observation to the rank of the intersection matrix of the irreducible components $Z_i$
of the complement $Z$ of $X$ in a smooth compactification.
\end{Exo}

%%% Local Variables:
%%% mode: latex
%%% TeX-master: "head"
%%% End:

\bigskip
%%%%%%%%%%%%%%%%%%%%%%%%%%%%%%%%%%%%%%%%%%%%%%%%%%%%%%%%%%%%%%%%%%%%%%%
%
%  Bibliography
%
%%%%%%%%%%%%%%%%%%%%%%%%%%%%%%%%%%%%%%%%%%%%%%%%%%%%%%%%%%%%%%%%%%%%%%%

\end{document}